\documentclass[11pt]{report}
\usepackage{graphics}
\usepackage{cancel}
\usepackage{amsmath}
\usepackage{graphicx}
\usepackage{epsfig}
\usepackage{graphicx}

\usepackage{amssymb}
\usepackage{amsmath, amsthm}

\newcommand{\beq}{\begin{equation}}\newcommand{\enq}{\end{equation}}
\newcommand{\C}{\mathbb{C}}\newcommand{\R}{\mathbb{R}}\newcommand{\im}{\text{Im }}
\newcommand{\N}{\mathbb{N}}\newcommand{\Z}{\mathbb{Z}}\newcommand{\eps}{\varepsilon}

\newcommand{\Raw}{\Rightarrow}
\newcommand{\raw}{\rightarrow}\newcommand{\new}{\newline}\newcommand{\bld}{\mathbf}
\newcommand{\tr}{\text{Tr }}\newcommand{\Pm}{\mathbb{P}}

\newcommand{\cone}{\begin{center}}\newcommand{\ctwo}{\end{center}}\newcommand{\stk}{\nu^{\star k}}

\newcommand{\bs}{\backslash}\newcommand{\beqa}{\begin{eqnarray*}}\newcommand{\enqa}{\end{eqnarray*}}
\newcommand{\half}{\frac{1}{2}}

\newcommand{\la}{\langle}\newcommand{\ra}{\rangle}

 \newcommand{\wxn}{\addtocontents{toc}{\protect\setcounter{tocdepth}{1}}}
\newcommand{\wxo}{\addtocontents{toc}{\protect\setcounter{tocdepth}{2}}}

\begin{document}

\cone
    \Large{The Cut-Off Phenomenon in Random Walks on Finite Groups}
\ctwo

\bigskip

\bigskip

\cone
 A thesis submitted to the National University of Ireland, Cork for the degree of Master of Science
\ctwo
    \author{Jeremiah Patrick McCarthy}

\cone
Supervisor: Dr. Stephen Wills
\ctwo

\cone
Head of Department: Prof. Martin Stynes
\ctwo

\bigskip
\cone
Department of Mathematics
\ctwo
\cone
College of Science, Engineering and Food Science
\ctwo
\cone
National University of Ireland, Cork
\ctwo

\bigskip

\cone
September 2010
\ctwo


\tableofcontents

\newpage
\section*{Abstract}
How many shuffles are needed to mix up a deck of cards? This question may be answered in the language of a random walk on the symmetric group, $S_{52}$. This generalises neatly to the study of \emph{random walks on finite groups} --- themselves a special class of Markov chains. \emph{Ergodic} random walks exhibit nice limiting behaviour, and both the quantitative and qualitative aspects of the convergence to this limiting behaviour is examined. A particular qualitative behaviour --- the \emph{cut-off phenomenon} --- occurs in many examples. For random walks exhibiting this behaviour, after a period of time,  convergence to the limiting behaviour is abrupt.

\bigskip

The aim of this thesis is to present the general theory of random walks on finite groups, with a particular emphasis on the cut-off phenomenon. It is an open problem to determine which random walks exhibit the cut-off phenomenon. There are various formulations of the cut-off phenomenon; the original --- that of \emph{variation distance} cut-off --- is considered here. At present, progress is made on this problem in a case-by-case basis. There are general techniques for attacking a particular case --- and many of these are presented here --- but there are no truly universal results.

\bigskip

Throughout the thesis, examples are used to demonstrate the theory. The last chapter presents some new heuristics developed by the author in the course of his studies.
\newpage
    \section*{Acknowledgements}
In the first instance, I would like to sincerely thank my supervisor Stephen Wills for proposing such an interesting area to work in. I am grateful for the fact that he provided assistance whenever I needed it --- be it in terms of mathematics, day-to-day life in the School of Mathematical Science, or practical advice for the writing of this thesis. I am very much looking forward to working with him on my Ph.D work.

\bigskip

I am sincerely grateful to Teresa Buckley and her team in the School office; never once were my problems left unresolved. I would also like to thank Martin Stynes and his engineering students (also students of Stephen Wills). As their tutor, I had to go over some old material: this undoubtedly saved me from some embarrassing mistakes in my thesis!

\bigskip

I would like to mention my friends in the Mathematics Research lab for enhancing my life in the School. I would like to share my appreciation of my grandmother Roses' cooking: it was a regular crutch over the last year! Also my housemates Ballsie, Swarley and Aliss --- it was always easy to unwind in their company after a long day at the sums.


\chapter{Introduction}
The question, how many shuffles are required to mix up a deck of cards , does not appear to have an obvious mathematical answer. Before any kind of analysis can be done, the terms  \textit{deck of cards}, \textit{shuffle} and \textit{mixed up} need  a precise mathematical realisation.
\new Consider a fresh deck of cards; in the order,  $K\heartsuit,Q\heartsuit,\dots,A\heartsuit,K\spadesuit,\dots,A\clubsuit$. In this order, each card can be labeled $1,\dots,52$, and given any arrangement of the deck, a permutation $\sigma:\{1,\dots,52\}\raw \{1,\dots,52\}$ can encode the arrangement:
\beq
\left(\begin{array}{cccc}
1 & 2 & \cdots & 52
\\ \sigma(1) & \sigma(2) & \cdots & \sigma(52)
\end{array}\right)
\enq
In the language of group theory, the deck of cards may be modelled by $S_{52}$.
\new A shuffle, meanwhile, takes the deck, and, \emph{independently\footnote{in general (!), one doesn't shuffle while looking at the labels on the cards. To be technical, not all functions $S:S_{52}\raw S_{52}$ are considered shuffles. For example, the `shuffle' swapping the positions of $A\heartsuit$ and  $A\spadesuit$ is not a shuffle.} of the arrangement of the deck}, permutes the cards. For example, a perfect cut shuffle takes off the top half of the cards and places it under the bottom half of the deck is a shuffle. It is not hard to see that  a shuffle is a function $S:S_{52}\raw S_{52}$, whose action is by multiplication by some $\sigma_S\in S_{52}$; i.e. $S(\sigma)=\sigma_S\sigma$. Indeed the perfect-cut shuffle is realised by multiplication by $(1,27)(2,28)\cdots(26,52)$.
\new Now the question of when is a deck mixed-up needs to addressed. In the first instance, it is always assumed that the deck started in some known order; e.g. the one given above.  Secondly, when is a deck totally random?
\newpage If one is handed a deck of cards, face down, and if each possible order of the cards is equally possible then the deck is considered random. It should be clear from  group theory, that if any perfect shuffle is repeated, then the deck will never get random in this sense. If the deck is always shuffled by $\sigma_S$, then after $k$ shuffles the deck will be in the order $\sigma_S^k$. Hence, to get random, there has to be some randomisation in how the deck is shuffled. As an example of a suitable randomisation, pick two distinct cards at random\footnote{to be careful maybe two distinct card \emph{positions}, e.g. top card, second card, etc.}, and let the shuffle swap the positions of these two cards. Assuming now that after a number of shuffles, every arrangement of the deck is approximately equally likely, various notions of `how close' the deck is to random may be formulated, and a clear definition of mixed-up may be given.

\bigskip

Consider the \textit{riffle shuffle}: at each step the deck is cut into two packs which are then \emph{riffled} together. A model for such shuffles on $n$ rather than just 52 cards, due to Gilbert, Shannon and (independently) Reeds, was completely analysed in a remarkable paper by Bayer \& Diaconis \cite{BD}. In this paper, a phenomenon called the \textit{cut-off phenomenon} was proven to occur for the riffle shuffle. Namely, for $n$ large, the deck is far from random in a certain sense after less than $t_n=(3\log_2n)/2$ shuffles, but close to random after more than $t_n$ shuffles: the transition from order to random takes place at about $t_n$ steps and it makes sense to say it takes $t_n$ steps to mix-up the cards. For the case $n=52$, seven shuffles are necessary and sufficient to mix up the cards.

\bigskip

Random walks on finite groups generalise card shuffling by replacing the symmetric group by any finite group. This thesis aims to present the general theory of random walks on finite groups, with an emphasis on the cut-off phenomenon. In particular, care has been shown to take no liberties with assumptions, and all the `obvious' elements of the theory are revisited and questioned. For example, Theorem \ref{eth} is standard in the field but almost all references do not carry the non-trivial proof. The questioning of `obvious' facets of the theory allowed some new perspectives.

\newpage

In making the thesis modest, some interesting and often powerful aspects of the theory have been omitted. The aforementioned riffle shuffle was not studied --- neither was the familiar over-hand shuffle. In fact, in terms of the development of the subject, the riffle shuffle is a pathological example. Despite its apparent complexity, the shuffle has been more or less completely understood and analysed by Bayer \& Diaconis, albeit through some deeper mathematics than the subject usually requires.
\new The Diaconis-Fourier theory is an attractive machinery in the field that is presented here. However it is only applied in two Abelian examples: neither of which needed require the full theory anyway. Its greatest success has been in the analysis of the random transposition shuffle, a random walk on the symmetric group, however the representation theory of the symmetric group is not covered here. Diaconis \cite{PD} is an excellent reference. A great survey of techniques, including those not mentioned here is \cite{sc}.
\new There are a number of interesting generalisations of random walks on groups, such as to homogenous spaces and Gelfand pairs. These are not covered here: Ceccherini-Silberstein \emph{et al} \cite{cecc} is an excellent book and pursues these areas.

\bigskip

Despite these restrictions, a great  variety of mathematical techniques are used. Probability, measure theory, representation theory, functional analysis, geometry and, naturally, group theory is used throughout the thesis. The cut-off phenomenon is not just a theory for random walks on groups,  it occurs for some more general  Markov chains also. A breakthrough in the theory of random walks on groups will surely have an impact for the Markov chain community. In his introduction, Chen \cite{chen} discusses a few examples where the existence of a cut-off has a significant impact for applications.

\bigskip

This first chapter introduces the general discrete time Markov chain theory on a finite set. Random walks on groups are introduced as a special class of Markov chains and necessary and sufficient conditions for a random walk to `get random' are developed.
\new Chapter 2 discusses what it means for a random walk to be `close to random'. A number of measures of closeness to random are introduced. A distinguished distance, namely the \textit{variation distance}, is identified as the conventional measure of closeness to random in this study. An interpretation of variation distance by Switzer is shown to be correct here. Much of the spectral analysis of the stochastic operator is done in this chapter and this yields upper bounds on the distance to random --- many related to the eigenvalues of the associated \emph{stochastic operator}. Next techniques for finding  lower bounds on the distance to random are discussed. Finally, methods of procuring bounds for these eigenvalues via the geometry of the group are presented.
\new Chapter 3 develops the representation theory of finite groups. In conjunction with Fourier analysis for finite groups, this machinery, so well pioneered by Diaconis, is a powerful technique for generating bounds on the distance to random. Here the full, general, theory is developed. Two Abelian examples, the simple walk on the circle and the simple walk with loops on the $n$-Cube, are analysed.
\new Chapter 4 introduces the cut-off phenomenon and its formulation. In particular, it is seen that the phenomenon is defined with respect to the limiting behavior of a  family of random walks on groups, $\{G_n:n\in \N\}$, as the size of the group increases to infinity ($n\raw\infty$). There is a discussion of the present understanding of the cut-off phenomenon, and reasons for its existence are mentioned.
\new Chapter 5 presents some probabilistic methods for bounding the distance to random. These powerful methods --- strong uniform times and coupling --- are occasionally very transparent and help explain \emph{why} cut-offs occur.
\new Finally in Chapter 6 some new viewpoints and generalisations are presented. Although the motion of a particle in a random walk \emph{is} random (in general, after $k$ steps the position of the particle is unknown), its distribution after $k$ steps is deterministic. Thus the random walk has the structure of a dynamical system. Here an attempt is made to develop this further. Also the question of whether or not the invertibility of the stochastic operator  has implications for a random walk is addressed. A study of invertible stochastic operators is, as far as this author knows, non-existent in the literature. A few basic properties and questions are explored. Finally, a conjecture of the author, namely that if the stochastic operator is invertible, then the cut-off phenomenon will not be exhibited, is explored and disproved.

\newpage
\section{Markov Chain Theory}
Essentially, a Markov Chain is a construction of a mathematical model for a certain type of discrete motion of a particle in a space. The particle begins at some initial point and at certain times $t_1,t_2,\dots$ moves to another point in the space chosen `at random'. The probability that the particle moves to a certain point $y$ at a time $t$ is dependent only upon its position $x$ at the previous time. This is the Markov property.

\bigskip

To formulate, let $X$ be a finite set. Denote by  $M_p(X)$ the probability measures on $X$. Let $\delta^x$ be the element of $M_p(X)$ which puts a measure of 1 on $x$ (and zero elsewhere). These \textit{Dirac measures},  $\{\delta^x:x\in X\}$, are the canonical basis for $\R^{|G|}\supseteq M_p(X)$. A probability measure $\nu\in M_p(X)$ is \textit{strict} if $\nu(x)>0$, for all $x\in X$.   Denote by $F(X)$ the complex functions on $X$ and $L(V)$ the linear operators on a vector space $V$. The similarly defined \textit{Dirac functions},  $\{\delta_x:x\in X\}$, are the canonical basis for $F(X)$.  With respect to this basis $P\in L(F(X))$ has a matrix representation $[p(x,y)]_{xy}$. $P\in L(F(X))$ is a \textit{stochastic  operator} if:

\begin{itemize}
\item[(i)] $p(x,y)\geq 0,\,\,\forall x,y$
\item[(ii)]  $\sum_yp(x,y)=1$, $\forall x$ \hfill(row sum is unity)
\end{itemize}


%
%
%
%
%


Given $\nu\in M_p(X)$, a stochastic operators $P$ acts on $\nu$ as $\nu P(x)=\sum_y\nu(y)p(y,x)$. Stochastic operators are readily characterised without using matrix elements  as being $M_p(X)$-stable in the sense that $M_p(X)P\subset M_p(X)$ if and only if $P$ is a stochastic operator. It is an immediate consequence that if $P$ and $Q$ are stochastic, then so is $PQ$.
 %
%
%





\newpage

\wxn\subsection{Definition}\wxo
Let $X$ be a finite set and $\nu\in M_p(X)$, $P$ a stochastic operator on $X$, and $(Y,\mathcal{A},\mu)$ a probability space. A sequence $\{\xi_k\}_{k=0}^n$ of random variables $\xi_k:Y\raw X$ are a \textit{Markov Chain with initial distribution $\nu$ and stochastic operator $P$}, if
\begin{itemize}
\item[(i)] $\mu(\xi_0=x_0)=\nu(x_0)$.
\item[(ii)] $\mu(\xi_{k+1}=x_{k+1}\,|\,\xi_0=x_0,\dots,\xi_k=x_k)=p(x_k,x_{k+1})$,\new  assuming $\mu(\xi_0=x_0,\dots,\xi_k=x_k)>0$.
\end{itemize}

\bigskip


If $\nu=\delta^x$ in (i) the Markov chain is said to \textit{start deterministically} at $x$. Condition (ii) is the \textit{Markov property}.  Subsequent references to a Markov Chain $\xi$ refer to a Markov Chain $\left(\{\xi_k\}_{k=0}^n,P,\nu\right)$.

In terms of existence, given $\nu$ and $P$, let $$Y:=X^{n+1}=\underbrace{X\times X\times\cdots \times X}_{n+1\text{ copies}}$$
Define $\xi_k:Y\raw X$ by $\xi_k(x_0,\dots,x_n)=x_k$ and $$\mu(x_0,\dots,x_n)=\nu(x_0)p(x_0,x_1)\cdots p(x_{n-1},x_m).$$
Then $\mu\in M_p(Y)$, and $\xi$ is a Markov Chain for $\nu$ and $P$.
\wxn\subsection{Example: Two State Markov Chain}\wxo
Consider the set $X=\{1,2\}$ and $\nu\in M_p(X)$. Suppose the probability of going from 1 to 2 is $p$ and the probability of going from 2 to 1 is $q$. Then the two state Markov chain has stochastic operator
\beqa
P=\left(\begin{array}{cc}1-p & p\\[2ex] q & 1-q\end{array}\right)
\enqa
for $p,q\in\,[0,1]$.

\begin{figure}[h]\cone\epsfig{figure=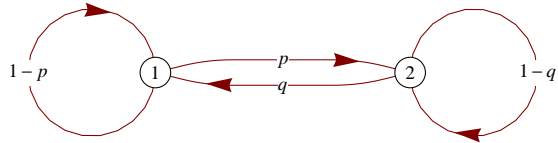}\ctwo\caption{A graphical representation of the two state Markov chain.}\end{figure}
\newpage
\section{Ergodic Theory\label{ergt}}
Ergodic theory is concerned with the longtime behaviour of a Markov
chain. A central question is for a given chain whether or not the
$\xi_k$ display limiting behaviour as $k\raw\infty$?  If
`$\xi_\infty$' exists, what is its distribution?

 \bigskip

 One possible debarring of the existence of a limit is periodicity. Consider a Markov chain $\xi$ on a set $X=X_0\cup X_1$ with $X_0\cap X_1=\emptyset$ and neither of the  $X_i=\emptyset$ for $i=1,2$. Suppose $\xi$ has the property that $\xi_{2k+i}\in X_i$, for $k\in\N_0$, $i=0,1$. Then `$\xi_\infty$' cannot exist in the obvious way. In a certain sense $\xi$ must be \textit{aperiodic} for limiting behaviour to exist.

 \bigskip

  Suppose $\xi$ is a Markov chain and the limit $\nu P^n\raw \theta$ exists. Loosely speaking, after a long time $N$, $\xi_N$ has distribution $\mu(\xi_N)\sim \theta$:
  \beqa
  \nu P^N\sim \theta
  \\\Raw \nu P^{N}P\sim\theta P
  \\\Raw \nu P^{N+1}\sim \theta P
  \enqa
  But $\nu P^{N+1}\sim \theta$ also and hence $\theta P\sim \theta$. So if `$\xi_\infty$' exists then its distribution $\theta$ may have the property $\theta P=\theta$. Such a distribution is said to be a \textit{stationary distribution for $P$}. Relaxing the supposition on `$\xi_\infty$' existing, do stationary distributions exist? Clearly they are left eigenvectors of eigenvalue 1 that have positive entries summing to 1.\newpage If $k(x)\in F(X)$ is any constant function then $Pk=k$ so $k$ is a right eigenfunction of eigenvalue 1. Let $u$ be a left eigenvector of eigenvalue 1. By the triangle inequality, $|u(x)|=|\sum_yu(y)p(y,x)|\leq\sum_y|u(y)|p(y,x)$. Now
  $$\sum_{z\in X}|u(z)|\leq \sum_{z\in X}\left(\sum_{y\in X}|u(y)|p(y,z)\right)=\sum_{y\in X}|u(y)|\underbrace{\left(\sum_{z\in X}p(y,z)\right)}_{=1}=\sum_{y\in X}|u(y)|$$
  Hence the inequality is an equality so $\sum_z\left(\sum_y|u(y)|p(y,z)-|u(z)|\right)=0$ is a sum of non-negative terms.  Hence $|u|P=|u|$, and by a scaling,  $\pi(x):=|u(x)|/\sum_y |u(y)|$, is a stationary distribution.


  \bigskip

  How many stationary distributions exist? Consider Markov Chains $\xi$ and $\zeta$ on disjoint finite sets $X$ and $Y$, with stochastic operators $P$ and $Q$. The block matrix
  \beq R=\left(\begin{array}{cc}P & 0\\0 & Q\end{array}\right)\enq
 is a stochastic operator on $X\cup Y$. If $\pi$ and $\theta$ are stationary distributions for $P$ and $Q$ then
 $$\phi_c=(c\pi,(1-c)\theta)\,,\,\,\,\,c\in[0,1]$$
 is an infinite family of stationary distributions for $R$. The dynamics of this walk are  that if the particle is in $X$ it stays in $X$, and vice versa for $Y$ (the graph of $R$ has two disconnected components). This example shows that, in general, the stationary distribution need not be unique. Rosenthal \cite{ros} shows that a sufficient condition for uniqueness is that the Markov chain $\xi$   has the property that every point is accessible from any other point; i.e. for all $\,x,y\in X$, there exists $r(x,y)\in\N$ such that $p^{(r(x,y))}(x,y)>0$. A Markov chain satisfying this property is said to be \textit{irreducible}.

\bigskip

So for the existence of a unique, stationary distribution it may be sufficient that the Markov chain is both aperiodic and irreducible. Call a stochastic operator $P$  \textit{ergodic} if there exists $n_0\in \N$ such that
$$p^{(n_0)}(x,y)>0\,,\,\,\forall x,y\in X$$
In fact, ergodicity is equivalent to aperiodic and irreducible  (see \cite{ros}\footnote{although aperiodic hasn't been defined here} Lemma 8.3.9), and the following theorem asserts that it is both a necessary and sufficient condition for the existence of a strict distribution for `$\xi_\infty$'. These precluding remarks  suggest the distribution of `$\xi_\infty$' is in fact stationary and unique, and indeed this will be seen to be the case. A nice, non-standard proof of this well-known theorem is to be found in \cite{cecc}.



\wxn\subsection{Markov Ergodic Theorem\label{erdt}}\wxo
\emph{A stochastic operator $P$ is ergodic if and only if there exists a strict $\pi\in M_p(X)$ such that
\beq\lim_{n\raw\infty} p^{(n)}(x,y)=\pi (y)\,,\,\,\forall x,\,y\in X\label{erg}\enq
In this case $\pi$ is the unique stationary distribution for $P$ $\bullet$}

\bigskip

In the special class of ergodic Markov chains, (\ref{erg}) indicates that statistically speaking, the system that evolves for a long time `forgets' its initial state. Another special class of Markov chains are \textit{reversible} Markov chains. A stochastic operator $P$ is \textit{reversible} if there exists a strict $\pi\in M_p(X)$ such that
\beq
\pi(x)p(x,y)=p(y,x)\pi(y)\,,\,\,\,\forall\,x,y\in X\label{dp}
\enq


This is equivalent to $D_\pi P=P^TD_\pi$ where $D_\pi$ is the diagonal matrix with $(x,x)$-component $\pi(x)$. Suppose further that $P$ is ergodic and (\ref{dp}) holds for some strict $\pi\in M_p(G)$. A quick calculation shows that then $\pi$ is the unique, strict, stationary distribution. The definition of a reversible chain  appears at odds with our interpretation of what reversible means. However, it may be shown (see \cite{cecc}) that the condition is equivalent to
\begin{itemize}
\item[(i)] $p(x,y)>0\Raw p(y,x)>0$
\item[(ii)] for all $n\in\N\,,\,\,x_0,x_1,\dots,x_n\in X$,
$$p(x_0,x_1)p(x_1,x_2)\cdots p(x_{n-1},x_n)p(x_n,x_0)=p(x_0,x_n)p(x_n,x_{n-1})\cdots p(x_1,x_0)$$
\end{itemize}


\begin{figure}[h]\cone\epsfig{figure=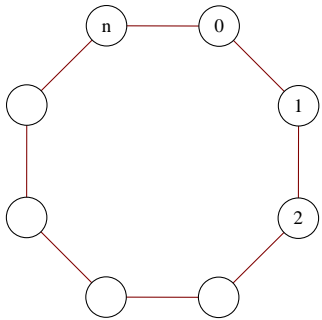}\ctwo\caption{For a reversible Markov Chain, the probability of going in a cycle from $0\raw 0$ is equal for clockwise and anti-clockwise orientations.}\end{figure}

\newpage
\section{Random Walks on Finite Groups}\wxn\subsection{Introduction\label{rwfgi}}\wxo A
particularly nice class of Markov chain is that of a \textit{random
walk on a group}. The particle moves from group element to group
element by choosing an element $h$ of the group `at random' and
moving to the product of $h$ and the present position $g$, i.e. the particle moves from $g$ to $hg$. To avoid trivialities, the random walk on the trivial group is not considered.
Naturally the group structure of the walk induces strong symmetry
conditions: this allows the generation of much stronger results than
that of general Markov chain theory.



\bigskip

To formulate, let $G$ be a finite group  of order $|G|$ and identity
$e$. Let $\nu\in M_p(G)$ and $(Y,\mu)$ be a probability space. Let
$\{\zeta_k\}_{k=0}^n: (Y,\mu)\raw G$ be a sequence of i.i.d.
random variables with distributions $\mu(\zeta_0=g_0)=\delta^e(g_0)$
and $\mu(\zeta_k=g)=\nu(g)$. The sequence of random variables
$\{\xi_k\}_{k=0}^n:(Y,\mu)\raw G$ \beq
\xi_k=\zeta_k\zeta_{k-1}\cdots\zeta_1\zeta_0 \enq is a
\textit{right-invariant random walk on $G$}.

\newpage

This construction makes $\xi$ into a Markov Chain on $G$ with
initial distribution $\delta^e$ and stochastic operator $P=p(s,t)$
is induced by the \textit{driving probability}, $\nu$: $p(s,t)=\nu(ts^{-1})$. The random walk is called
\emph{right invariant} because $p(s,t)=p(sh,th)$. This is obvious as
\beqa
p(sh,th)=\nu(th(sh)^{-1})=\nu(ts^{-1})=p(s,t)
\enqa
\subsubsection{Example: Card Shuffling}
Card shuffling provides the motivation  for the study of random walks on groups and remains the canonical example. Everyday shuffles such as the overhand shuffle or the riffle shuffle, as well as simpler but more tractable examples such as top-to-random or random transpositions all have the structure of a random walk on $S_{52}$. Each shuffle may be realised as sampling from a probability distribution $\nu\in M_p(S_{52})$. For example, consider the case of repeated random transpositions. A \emph{random transposition} consists chooses two cards at random (with replacement) from the deck and swapping the positions of these two cards. Suppose without loss of generality that the first card chosen is the ace of spades. The probability of choosing the ace of spaces again is 1/52. Swapping the ace the spades with itself leaves the deck unchanged. The choice of the first card is independent hence the probability that the shuffle leaves the deck unchanged is 1/52. What is the probability of transposing two given (distinct) cards?  Consider, again without loss of generality, the probability of transposing  the ace of spades and the ace of hearts. There are two ways this may be achieved: choose $A\spadesuit$-$A\heartsuit$ or choose $A\heartsuit$-$A\spadesuit$. Both of these have probability of 1/52$^2$. Any other given shuffle (not leaving the deck unchanged or transposing two cards) is impossible. Hence the shuffle may be modelled as sampling by
\beqa
\nu(s):=\left\{\begin{array}{cc}
1/52 & \text{ if }s=e
\\ 2/52^2 & \text{ if $s$ is a transposition }
\\ 0 & \text{ otherwise }
\end{array}\right.
\enqa

\newpage

  It is a straightforward calculation to show that the stochastic operator of a random walk on a group  is
doubly stochastic --- column sums are also 1. As a corollary, the
uniform distribution, $\pi(g)=1/|G|$, is a strict, stationary
distribution. To keep terminology to a minimum, the uniform distribution shall be referred to as the \textit{random distribution} and conversely $\pi$ will refer to this random distribution.



%
If $\Sigma=\text{supp }(\nu)$, then, in general, $\xi_k\in \Sigma^k$ however  if $\la\Sigma\ra=G$ and $e\in\Sigma$ then certainly $\Sigma^k\subset \Sigma^l$, for any $k\leq l$. Indeed:
$$
\{e\}=\Sigma^0\subset \Sigma\subset \Sigma^2\subset \cdots\subset \Sigma^T=G
$$
where $T$ is called the \textit{cover time} of the walk. In this case $P$ is ergodic with $n_0= T$. From Section \ref{ergt}, it is known
that `$\xi_\infty$' exists in a nice way if the stochastic operator
$P$ is ergodic.    Conveniently, this condition may be translated into
a condition on the driving probability on the group, $\nu$. The below theorem falls under the category of a `folklore theorem' in that almost all references  refer to the proof in older hard-to-source references --- if at all. A proof outline is given by Fountoulakis \cite{bhat} in his lecture notes but here a full proof is given.

\wxn\subsection{Ergodic Theorem for Random Walks on Groups\label{eth}}\wxo \emph{Let $G$ be a group and
$\nu\in M_p(G)$ with support $\Sigma$. A right-invariant random walk
on $G$ is ergodic if and only if $\Sigma\not\subset K$ for any proper subgroup $K$  of $G$
and $\Sigma\not\subset Hx$ for any coset of any proper normal
subgroup $H\lhd G$.
\new In this case, $\pi $ is the unique, strict stationary distribution for
$P$.}

\begin{proof}
Assume $\Sigma\subset K$ a proper subgroup of $G$.
$\la\Sigma\ra\subset K$ by closure in $K$; hence $\xi_k\in K$,
for all $k\in \N$. Let, $s\in K$, $t\not\in K$. Now for all  $n\in \N$,
$p^{(n)}(s,t)=0$. Hence $P$ is \emph{not} ergodic.
\new Assume $\Sigma \subset Hx$ for some coset of a proper normal subgroup $H\lhd G$. Now $\xi_0\in He$ and $\xi_1\in HxHe=Hx$, so by induction $\xi_n\in (Hx)^n=Hx^n$, for all $n\in
\N$.   Let $n\in\N$. Let $s\in G\bs Hx^n$: $p^{(n)}(e,s)=0$. Hence $P$
is \emph{not} ergodic.

\bigskip

Assume now $\Sigma\not\subset K$ a proper subgroup of $G$ and
$\Sigma \not\subset Hx$ for any coset of any proper normal subgroup
$H\lhd G$.
\new Clearly the inclusions $\Sigma \subset \la\Sigma \ra\subset G$ hold with $\la \Sigma\ra$  a subgroup of $G$.  By assumption $\Sigma$ does not lie in a proper subgroup hence $\la\Sigma\ra=G$. Hence for all $s,t\in G$, there exists $n(s,t)\in \N$
such that $p^{(n(s,t))}(s,t)>0$.
\new Let $L_\Sigma(e):=\{(\sigma_{i_1},\dots,\sigma_{i_N}):e=\sigma_{i_1}\cdots\sigma_{i_N}\,;\sigma_{i_m}\dots\sigma_{i_n}\neq
e\,,n-m< N-1\,;\sigma_{i_j}\in\Sigma\}$ be the set of
all distinct minimal $\Sigma$-presentations of $e$.


\bigskip

\textbf{Claim 1: } If $|L_\Sigma(e)|=1$, $G=\Z_{|G|}$ and $\Sigma$
is in a coset of a proper normal subgroup.
\new\emph{Proof}. If $|L_\Sigma(e)|=1$ there is only one minimal $\Sigma$-presentations of $e$. But
  $\sigma_1^{o(\sigma_1)}$ and $\sigma_2^{o(\sigma_2)}$ are
two distinct minimal $\Sigma$-presentations of $e$. Hence
$\sigma_1=\sigma_2$. Hence $\Sigma=\{\sigma\}$. But
$\la\Sigma\ra=G$, hence $G$ is cyclic and in particular
$\Sigma\subset \{e\}\sigma$ the coset of the proper normal subgroup
$\{e\}\,\,\,\bullet$



\textbf{Claim 2:} Assume $|L_\Sigma(e)|>1$. If $\Sigma$ is not
contained in a coset of a proper normal subgroup of $G$, then, where $L$ is the set of word lengths of the elements of $L_\Sigma(e)$, $\gcd
L=1$. \new \emph{Proof.} Suppose $\gcd L=k>1$.
Then every $\Sigma$-presentation of $e$ has length $0$ mod $k$. Let
$N_k \subset G$ be the subgroup generated by all elements of $G$
with a length $0$ mod $k$  $\Sigma$-presentation. Clearly $e\in
N_k$. Let $t\in G$. Suppose $t$ has a  length $p$ mod $k$
$\Sigma$-presentation. Then $t^{-1}$ has a length $-p$ mod $k$
$\Sigma$-presentation since $t^{-1}t=e$ has length
$0$ mod $k$ $\Sigma$-presentation. Let $n\in N_k$. By
definition, $n$ has a length $0$ mod $k$  $\Sigma$-presentation and
so $t^{-1}nt$ has a length $0$ mod $k$  $\Sigma$-presentation. So
$N_k$ is normal.
\new Let $\sigma\in\Sigma$ and suppose
$\sigma\in N_k$. Then
  \beqa
   \sigma\sigma^{-1}=e=(\sigma_{i_1}\cdots \sigma_{i_{qk}})(\sigma_{j_1}\cdots \sigma_{j_{lk-1}})
  \enqa
   that is $e$ would have a length $-1$ mod $k$ $\Sigma$-presentation, which is not allowed. Hence $\sigma\not\in N_k$, so $N_k$ is a proper normal subgroup of $G$.  \new Let $\sigma_1\in\Sigma$. Then $\sigma\sigma^{-1}_1\in N_k$ for all $\sigma\in \Sigma$ as $\Sigma$-presentations of any $\sigma^{-1}$ have length $-1$ mod $k$. Hence $\Sigma \subset N_k\sigma_1$  and this contradicts the assumption on $\Sigma$. Hence $\gcd
L_\Sigma(e)=1\,\,\,\bullet$



 Let $S$  be the set of lengths of all\footnote{not just minimal presentations} distinct
$\Sigma$-presentations of $e$.  As $L\subset S$, $\gcd
S=1$. Hence there exist $l_1,\dots,l_m\in S$, $k_i\in\Z$ such that \cite{H}: \beq k_1l_1+\cdots+k_ml_m=1 \enq Let $l\in S$ and $n(e,s)$
as above.
\new Let
\beq M=l_1|k_1|+\cdots+l_m|k_m| \enq and \beq n_0(e,s)=lM+n(e,s)
\enq If $n\geq n_0(e,s)$, and letting \begin{align*}
r&=\left\lfloor\frac{n-n(e,s)}{l}\right\rfloor\text{ , and } \\
n&=n(e,s)+rl+a \end{align*} where $0\leq a<l$ and $r\geq M$. Now as
\beqa\sum_{i=1}^m k_il_i=1 \text{ , and }\sum_{i=1}^ml_i|k_i|=M,
\enqa
$n$ may be written
\begin{align*}
n & =n(e,s)+rl\underbrace{-lM+l\left(\sum_{i=1}^ml_i|k_i|\right)}_{=0}+a\left(\sum_{i=1}^m k_il_i\right)
\\ & =(r-M)l+\sum_{i=1}^m(l|k_i|+ak_i)l_i+n(e,s)
\end{align*}
where the $(l|k_i|+ak_i)\geq0$. Let $x,y,\lambda\in\N$. Note that the probability of going from $s$ to $t$ in $x+\lambda y$ steps is certainly greater than going from $s$ to $t$ in $x$ steps and returning to $t$ every $y$ steps $\lambda$ times:
 \beq
 p^{(x+\lambda y)}(s,t)\geq p^{(x)}(s,t)\left(p^{(y)}(t,t)\right)^\lambda
 \enq

 Hence as $l,l_i\in S$ (so that $p^{(l)}(e,e)>0$) and $p^{(n(e,s))}(e,s)>0$;
\beqa p^{(n)}(e,s)\geq
\left(p^{(l)}(e,e)\right)^{r-M}\left[\prod_{i=1}^m\left(p^{(l_i)}(e,e)\right)^{l|k_i|+ak_i}\right]p^{(n(e,s))}(e,s)>0
\enqa Now let $n_0$ be the maximum of $n_0(e,s)$ as $s$ runs over $G$. Let $s,t\in G$. By right
invariance \beqa
p^{(n)}(s,t)=p^{(n)}(e,ts^{-1})>0\text{ , for
}n>n_0 \enqa Hence $P$ is ergodic $\bullet$
\end{proof}


%


\chapter{Distance to Random\label{dist}}
\section{Introduction}
The previous chapter demonstrates that under mild conditions a random walk on a group converges to the random distribution. Therefore, initially the walk is `far' from random and eventually the walk is `close' to random. An appropriate question therefore, is given a control $\eps>0$,
how large should $k$ be so that the walk is $\eps$-close to random after $k$ steps? The first problem here is to have a measure of `close to random'. This chapter introduces a few measures of `closeness to random', discusses the relationship between them and presents some bounds. In the rest of the work, all walks are assumed ergodic unless stated otherwise.

\bigskip

Let $\nu$ and $\mu\in M_p(G)$. The
\textit{convolution} of $\nu$ and $\mu$ is the probability \beq
\nu\star \mu(s):=\sum_{t\in G}\nu(st^{-1})\mu(t). \enq In particular
denote $\nu^{\star n+1}:=\nu\star \nu^{\star
n}$. The distribution of a random walk after one step is given by $\nu$. If $s\in G$, then the walk can go to $s$ in two steps by going to \emph{some} $t\in G$ after one step and going from there to $s$ in the next. The probability of going from $t$ to $s$ is given by the probability of choosing $st^{-1}$, i.e. $\nu(st^{-1})$. By summing over all intermediate steps $t\in G$, and noting that $\nu\star\delta^e=\nu$, it is seen that if $\{\xi_k\}_{k=0}^n$ is a
random walk on $G$ driven by $\nu$, then $\nu^{\star k}$ is the
probability distribution of $\xi_k$. In terms of the stochastic operator induced by $\nu\in M_p(G)$ , $P$, given any $\mu\in M_p(G)$, $\mu P=\nu\star \mu$.


\section{Measures of Randomness}

The preceding remarks indicate that $\nu^{\star k}\raw \pi$ thus a measure of closeness to random can be defined by defining a metric on $M_p(G)$ or putting a norm on $\R^{|G|}\supseteq M_p(G)$. Then a precise mathematical question may be asked: given $\eps>0$, how large should $k$ be so that $\|\nu^{\star k}-\pi\|<\eps$ or $d(\nu^{\star k},\pi)<\eps$? Straightaway it is clear that any of the $p$-norms may be used. Also multiples of $p$-norms may be used, for example, Diaconis \& Saloff-Coste \cite{chensixteen} introduce the distance $d_p(k):=|G|^{1-1/p}\|\nu^{\star k}-\pi\|_p$.

 \bigskip

 Another notion of closeness to random, although not a metric, is that of \textit{separation distance}:
\beq
s(k):=|G|\max_{t\in G}\left\{\frac{1}{|G|}-\nu^{\star k}(t)\right\}
\enq
Clearly $s(k)\in[0,1]$ with $s(k)=1$ if and only if $\nu^{\star k}(g)=0$ for some $g$; and $s(k)=0$ if and only if $\nu^{\star k}=\pi$. The separation distance is submultiplicative in the sense that $s(k+l)\leq s(k)s(l)$, for $k,l\in\N$ \cite{chenthree}. This immediately implies that $s(nk)\leq [s(k)]^n$. Suppose however that $\nu^{\star k}(g)=0$ for some $g\in G$. Then $s(k)=1$ and $s(nk)\leq 1$ which is useless. However because the walk is ergodic there exists a time $n_0$ when $\nu^{\star k}$ is supported on the entire group. Let $L:=\min\{\nu^{\star n_o}(s):s\in G\}$.  Then $s(n_0)=(1-|G|L)$, thence $s(kn_0)\leq(1-|G|L)^k$. An example where this bound is easily applied is the simple walk on $\Z_n$, $n$ odd, where $\nu(\pm1)=1/2$. Then $n_0=n-1$, $L=2^{1-n}$ and thence $s(k(n-1))\leq(1-n.2^{1-n})^k$.

\bigskip

A further measure of randomness is that of the average Shannon Entropy of the distribution; $H(\mu)=\sum_{t}\mu(t)\log\left(1/\mu(t)\right)$. A quick calculation shows that $H(\delta^e)=0$, $H(\pi)=\log|G|$; and also that $H(\nu^{\star k})$ increases to $\log |G|$ monotonically \cite{info}. Therefore $\sigma(k):=\log |G|-H(\nu^{\star k})$ is a measure of closeness to random.  A lower bound, adapted from \cite{chenone}, is $\sigma(k)\geq (1-k)\log|G|+k\sigma(1)$.

 \bigskip

 The default measure of closeness to random in this work, however, is \textit{variation distance}. If $\mu,\nu\in M_p(G)$, their \textit{variation distance} is
\beq
\|\mu-\nu\|:=\max_{A\subset G}|\mu(A)-\nu(A)|
\enq
\newpage
Diaconis \cite{PD} notes  an interpretation of variation distance of Paul Switzer. Consider $\mu$, $\nu\in M_p(G)$. Given a single observation of $G$, sampled from $\mu$ or $\nu$ with probability $1/2$, guess whether the observation, $o$, was sampled from $\mu$ or $\nu$. The classical strategy presented here gives the probability of being correct as $1/2(1+\|\mu-\nu\|)$:
\begin{enumerate}
\item Evaluate $\mu(o)$ and $\nu(o)$.
\item If $\mu(o)\geq\nu(o)$, choose $\mu$.
\item If $\nu(o)>\mu(o)$, choose $\nu$.
\end{enumerate}
To see this is true, let $\{\mu>\nu\}$  be the set $\{t\in G:\mu(t)>\nu(t)\}$.
Suppose $o$ is sampled from $\mu$. Then the strategy is correct if $o\in\{\mu=\nu\}$ or $o\in\{\mu>\nu\}$:
$$\Pm[\text{guessing correctly}\,|\,\mu]=\Pm[o\in\{\mu=\nu\}\,|\,\mu]+\Pm[o\in\{\mu>\nu\}\,|\,\mu]$$
with a similar expression for $\Pm[\text{guessing correctly}\,|\,\nu]$.
Note that $\Pm[o\in\{\mu=\nu\}]=\mu(\{\mu=\nu\})=\nu(\{\mu=\nu\})$ and also $\Pm[o\in\{\mu>\nu\}\,|\,\mu]=\mu(\{\mu>\nu\})$ (and similar for $o\in\{\mu<\nu\}$). Thus
\begin{align*}
\Pm[\text{guessing correctly}] &=\half\Pm[\text{guessing correctly}\,|\,\mu]+\half\Pm[\text{guessing correctly}\,|\,\nu]
\\&=\half\left(\nu(\{\mu=\nu\})+\mu(\{\mu>\nu\})\right)+\half\left(\nu(\{\mu<\nu\})\right)
\end{align*}
It is easily shown that
\beqa\|\mu-\nu\|=\mu\left(\{\mu>\nu\}\right)-\nu\left(\{\mu>\nu\}\right).
\enqa
Hence
$$
\Pm[\text{guessing correctly}]=\half\left(\underbrace{\nu(\{\mu=\nu\})+\nu(\{\mu>\nu\})+\nu(\{\mu<\nu\})}_{=1}+\|\mu-\nu\|)\right).$$

Also the separation distance controls the variation distance as
\beqa
\|\nu^{\star k}-\pi\|=\sum_{t\in\{\nu^{\star k}<\pi \}}\left(\frac{1}{|G|}-\nu^{\star k}(t)\right)\leq s(k).
\enqa
\newpage
It is a straightforward exercise, however, to show that $\|\mu-\nu\|$ is simply half of the usual $l^1$-distance $\|\mu-\nu\|_1$. Hence, with $P$ doubly stochastic ($\|P\|_{l^1\raw l^1}=1$ as column sums are 1), the quick calculation
\beqa
\|\nu^{\star k+1}-\pi\|_1=\|(\nu^{\star k}-\pi)P\|_1\leq \|\nu^{\star k}-\pi\|_1\|P\|_{l^1\raw l^1}=\|\nu^{\star k}-\pi\|_1
\enqa
shows that $\|\nu^{\star k}-\pi\|$ is decreasing in $k$.

\bigskip

At this juncture Aldous \cite{chenone} denotes by $\tau(\eps)$ the time to get $\eps$-close to random: $\min \{k:\|\nu^{\star k}-\pi\|<\eps\}$. Call $\tau:=\tau(1/2e)$ the \textit{mixing time}.   The reason the random walk driven by $\nu\in M_p(G)$ is defined to
start deterministically at $e$ is because due to right-invariance a
random walk driven by the same measure starting deterministically at
$g\neq e$ will converge to random at the same rate. Also, if $\xi_0$ is distributed as $\theta=\sum_t a_t\delta^t$, then the walk looks like $\bigoplus_t a_t \xi^t$ where   $\xi^t$ is the walk which begins deterministically at $t$. All these constituent  walks converge at the same rate, however, as might be expected:
\begin{align}
\|\theta P^k-\pi\| &=\half \sum_{s\in G} \left|\left(\sum_{t\in G}a_t \delta^t P^k(s)\right)-\pi(s)\right|=\half \sum_{s\in G}\left|\sum_{t\in G} a_t\left(\delta^t P^k(s)-\pi(s)\right)\right|\nonumber
\\&\leq \half \sum_{s\in G}\sum_{t\in G}a_t|\delta^tP^k(s)-\pi(s)|=\sum_{t\in G}a_t\left(\half\sum_{s\in G}|\delta^tP^k(s)-\pi(s)|\right)\nonumber
\\ & \leq \|\nu^{\star k}-\pi\|\label{slow}
\end{align}
Certainly there is equality if $\theta$ is a Dirac measure or the random distribution, $\pi$.

\section{Spectral Analysis\label{spectral}}
In the case of reversible random walks, where $\pi$ is the random distribution, $\pi(g)p(g,h)=p(h,g)\pi(h)$. Hence the driving probability is symmetric:
\beqa
p(g,h)=p(h,g)\Leftrightarrow \nu(hg^{-1})=\nu(gh^{-1})\Leftrightarrow \nu(s)=\nu(s^{-1})\,,\,\,\forall\,s\in G
\enqa
Also in the $\{\delta_t:t\in G\}$ basis the matrix representation of the stochastic operator is symmetric: $p(x,y)=p(y,x)$. Let $(\,|\,)$ be the  inner product on $F(G)$:
$$(\phi|\psi):=\frac{1}{|G|}\sum_{s\in G}\phi(s)\psi(s)^\star$$
When the walk is reversible:
\begin{align*}
(P\phi|\psi)&=\frac{1}{|G|}\sum_{s\in G}\left(\sum_{t\in G}p(s,t)\phi(t)\right)\psi(t)^\star
\\ & = \frac{1}{|G|}\sum_{t\in G}\phi(t)\left(\sum_{s\in G}p(t,s)\psi(s)\right)^\star=(\phi|P \psi),
\end{align*}
and so the stochastic operator is self-adjoint. By the spectral theorem for self-adjoint maps $P$ has an (left) eigenbasis $\mathcal{B}=\{u_1,\dots,u_{|G|}\}$. Suppose further that $\mathcal{B}$ is normalised such that $\delta^e=\sum a_tu_t$ with $u_1=\pi$ and $a_1=1$ (in fact for any $\theta\in M_p(G)$ this normalisation is unique. Let $v\in\R^n$. Call the sum of the entries of $v$ its \textit{weight}. The eigenvectors $u_t$, $t\neq 1$, are orthogonal to $\pi$.  Thence these eigenvectors have weight 0 so in order for the linear combination to be a probability distribution the weight needs to be 1, hence $a_1$ must be 1.). If $P$ is ergodic, then the eigenvalue $1$ has multiplicity 1. A quick calculation shows that if $\lambda_1=1$, then also $|\lambda_t|\leq1$, for all $t\neq 1$. Using an elegant graph-theoretic argument, Ceccherini-Silberstein \emph{et al} \cite{cecc} show that if $P$ is ergodic then $-1$ is not an eigenvalue. Therefore in the case of reversible walks (real eigenvalues), $|\lambda_t|<1$, for all $t\neq 1$ (this is also a consequence of the Perron-Frobenius Theorem), and then
 \beq
\nu^{\star k}=\delta^e P^k=\pi+\sum_{t\neq 1}a_t\lambda^k_tu_t
 \enq
Therefore, letting $ \lambda_{\star}:=\max\{|\lambda_t|:t\neq 1\}$;
\begin{align*}
\|\nu^{\star k}-\pi\|&=\half\left|\sum_{t\neq1}a_t\lambda_t^ku_t\right|=\half\sum_{s\in G}\left|\sum_{t\neq 1}a_t\lambda_t^ku_t(s)\right|
\\ & \leq \half\sum_{s\in G}\sum_{t\neq 1}|a_t||\lambda_t|^k|u_t(s)|
\\ & \leq \lambda_{\star}^k\underbrace{\half\sum_{s\in G}\sum_{t\neq 1}|a_t||u_t(s)|}_{=C}=C\lambda_\star^k
\end{align*}
Hence the rate of convergence is controlled by the second highest eigenvalue in magnitude. In Corollary \ref{she} an explicit $C$ is given. The importance of the second largest eigenvalue is a mantra in Markov chain theory, however it is only in the reversible case that the importance is so obvious.\newpage Suppose now that $P$ is a not-necessarily-reversible stochastic operator. Following \cite{2EV}, put $P$   in Jordan normal form:
\beqa
P=\left(\begin{array}{cccc}1 & & & 0\\ & J_2 & & \\ & & \ddots & \\ 0 & & & J_m\end{array}\right)
\enqa
where the Jordan blocks $J_i$ have form:
\beqa
J_i=\left(\begin{array}{cccc}\lambda_i & 1 & & 0\\ 0 & \lambda_i & \ddots & \\ &\ddots & \ddots & 1 \\ 0 & \cdots & 0 & \lambda_i\end{array}\right)
\enqa
and have size equal to the algebraic multiplicity of $\lambda_i$. Note the first entry of $P$ will be just 1 as 1 is an eigenvalue of multiplicity 1. The Jordan block $J_i$ is  the sum of the diagonal matrix $\lambda_iI$ and the superdiagonal, and thus nilpotent, matrix $N_i$.
   With $P^n=\text{diag}(1,J_1^n,\dots,J_m^n)$, and noting $N_i^{d_i}=0$ where $d_i$ is the multiplicity of $\lambda_i$;
\begin{align*}
J_i^k &=(\lambda_i I+N_i)^k=\sum_{j=0}^k{k\choose j}\lambda_i^{k-j}N_i^j=\sum_{j=0}^{d_i-1}{k\choose j}\lambda_i^{k-j}N_i^j.
\end{align*}
Now, for $j< d_i$, $N_i^j$ is the matrix with ones on the  $j$th diagonal above the main diagonal. Hence $J_i^k$ is a matrix whose lower diagonal entries are zero and have equal entries along this `$j$th diagonal', namely
\beqa
(J_i^k)_j={k\choose j}\lambda_i^{k-j}
\enqa
Hence the magnitude of the entries along the $j$th diagonal is bounded by (as $|\lambda_i|<1$):
\beqa
\left|(J_i^k)_j\right|\leq |\lambda_i|^k{k\choose j}
\enqa
The remaining manipulations are dependent on the relation of $k$ to $d_i$. Assuming $k>2d_i$ for example:
\beqa
\left|(J_i^k)_j\right|\leq |\lambda_i|^k{k\choose d_i}
\enqa
In Jordan normal form, $P$ converges to the matrix with $1$ in the $(1,1)$ entry and zero elsewhere. Clearly it is the block corresponding to the second largest eigenvalue in magnitude which is the slowest to converge and hence this eigenvalue controls convergence.

\bigskip

Taking the approach of \cite{cecc}, more explicit bounds for the reversible case may be found.
If the walk is reversible then $P$ has an (right) orthonormal basis $\mathcal{B}=\{v_t:t\in G\}$ with corresponding eigenvalues $\{\lambda_t:t\in G\}$. Let $v_1$ be the constant function with value 1 (so that $\lambda_1=1$). Put $\Lambda=\text{diag}(\lambda_1,\dots,\lambda_{|G|})$. Now
\beqa
Pv_s(g)=\sum_tp(g,t)v_s(t)=v_s(g)\lambda_s\Leftrightarrow PU=U\Lambda
\enqa
where $U=[v_1|\cdots|v_{|G|}]$. From orthonormality
\beqa
(v_s|v_h)=\frac{1}{|G|}\sum_tv_s(t)v_h(t)=\delta_s(h)\Leftrightarrow U^TU=|G|I
\enqa
 As a matrix of eigenvectors, $U$ is invertible with $U^{-1}=U^T/|G|$. Hence $P=U\Lambda U^T/|G|$, and so:
\beqa
P^k=\frac{1}{|G|^k}\underbrace{U\Lambda \underbrace{U^TU}_{=|G|}\Lambda U^T\cdots \Lambda U^T}_{k\text{ copies}}=U \Lambda^k U^T/|G|
\enqa
Or, in terms of coordinates,
\beqa
 p^{(k)}(g,h)=\frac{1}{|G|}\sum_{t\in G}v_t(g)\lambda_t^k v_t(h)
\enqa


 \wxn\subsection{Proposition}\wxo
 \emph{Suppose $\nu$ is symmetric. Then in the notation above}
 \beq
 \|\nu^{\star k}-\pi\|^2_2=\frac{1}{|G|}\sum_{t\neq 1}\lambda_t^{2k}v_t(e)^2
 \enq
 \begin{proof}
By definition
 \begin{align*}
\|\nu^{\star k}-\pi\|^2_2&=\sum_{s\in G}(\nu^{\star k}(s)-\pi(s))^2
 \\  &=\sum_{s\in G}\left(\sum_{t\neq 1}v_t(e)\lambda_t^kv_t(s)/|G|\right)^2
 \\ &=\sum_{t_1,t_2\neq 1}v_{t_1}(e)v_{t_2}(e)\lambda_{t_1}^k\lambda_{t_2}^k\sum_{s\in G} v_{t_1}(s)v_{t_2}(s)/|G|^2
 \end{align*}
 But $U^TU/|G|=I$; equivalently
 \beqa
 \sum_{s\in G}v_{t_1}(s)v_{t_2}(s)/|G|=\delta_{t_1}(t_2)
\enqa
and so
\beqa \|\nu^{\star k}-\pi\|^2_2=\frac{1}{|G|}\sum_{t\neq 1}v_t(e)^2\lambda_t^{2k}\,\,\,\bullet
 \enqa
\end{proof}

 \wxn\subsection[Upper Bound Lemma]{Corollary: Upper Bound Lemma}\wxo
\emph{Using the same notation, where $\|\cdot\|$ is the variation distance:}
\beq
 \|\nu^{\star k}-\pi\|^2\leq\frac{1}{4}\sum_{t\neq 1}v_t(e)^2\lambda_t^{2k}
 \enq
 \begin{proof}
The proof is a rudimentary application of the Cauchy-Schwarz Inequality:
 \begin{align*}
 \|\nu^{\star k}-\pi\|^2&=\frac{1}{4}\|\nu^{\star k}-\pi\|^2_1
 \\&=\frac{1}{4}\left(\sum_{t\in G}\left\{|\nu^{\star k}(t)-\pi(t)|\sqrt{|G|}\right\}\cdot\frac{1}{\sqrt{|G|}}\right)^2
 \\&\leq\frac{1}{4}\underbrace{\left(\sum_{t\in G}|\nu^{\star k}(t)-\pi(t)|^2|G|\right)}_{=|G|\|\nu^{\star k}-\pi\|_2^2}\left(\sum_{t\in G}\frac{1}{|G|}\right)\,\,\,\bullet
 \end{align*}
\end{proof}

 \wxn\subsection{Corollary\label{she}}\wxo
\emph{In the same notation:}
 \beq
 \|\nu^{\star k}-\pi\|^2\leq \frac{|G|-1}{4}(\lambda_\star)^{2k}
 \enq

 \begin{proof}
Since $|\lambda_t|\leq \lambda_\star$ for all $t\neq1$,
 \beqa
 \|\nu^{\star k}-\pi\|^2\leq\frac{1}{4}\sum_{t\neq 1}(\lambda_t)^{2k}v_t(e)^2\leq\frac{(\lambda_\star)^{2k}}{4}\sum_{t\neq 1}v_t(e)^2
 \enqa
Note that the eigenvectors of symmetric matrices can be chosen to be real-valued \cite{AbMagMA}, so that $v_t\overline{v_t}=v_t^2$. Also $UU^T=|G|I$ and hence $U^TU=|G| I$ thus
 \begin{align*}
 \sum_{t\in G}v_t(e)^2&=|G|
\\ v_1(e)^2+\sum_{t\neq 1}v_t(e)^2&=1+\sum_{t\neq 1}v_t(e)^2=|G|
\,\,\,\bullet
 \end{align*}
\end{proof}

When $\nu$ is symmetric, the associated stochastic operator, $P$, is symmetric and hence has real eigenvalues which can be ordered $1=\lambda_1>\lambda_1\geq \cdots\geq \lambda_{|G|}>-1$. So now $\lambda_\star=|\lambda_2|$ or $|\lambda_{|G|}|$. Of course, if the spectrum of $P$ can be calculated then these bounds are immediately applicable, however more often one must do with estimates. Diaconis and Saloff-Coste \cite{chensixteen} has many examples. Lemma 1 in that paper is a standard result in the field and is proved by consideration of the probability $\nu^\prime=(\nu-\nu(e)\delta^e)/(1-\nu(e))$ however a quick application of Gershgorin's circle theorem \cite{Lewis} shows the $\lambda_{|G|}\geq-1+2\nu(e)$ result also. As the Gershgorin result is mentioned in the sequel, and not typically used by the random walk community, it is presented here:
\wxn\subsection{Gershgorin's Circle Theorem}\wxo
\emph{Let $A$ be a complex $n\times n$ matrix with entries $a_{ij}$. Let $R_i=\sum_{j\neq i}|a_{ij}|$ be the sum of the absolute values of the entries in the $i$th row, excluding the diagonal element. If $B[a_{ii},R_i]$ is the closed disc centered at $a_{ii}$ with radius $R_i$, then each of the eigenvalues of $A$ is contained in at least one of the $B[a_{ii},R_i]$.}
\begin{proof}
Let $\lambda$ be an eigenvalue of $A$ with eigenvector $v$. Let $|v(k)|=\max_j|v(j)|$. Now
\beqa
Av(k)=\sum_{j=1}^na_{kj}v(j)=\lambda v(k).
\enqa
That is
\beqa
\sum_{j\neq k}a_{kj}v(j)=\lambda v(k)-a_{kk}v(k).
\enqa
Divide both sides by $v(k)$;
\beqa
\lambda-a_{kk}=\frac{\sum_{j\neq k}a_{kj}v(j)}{v(k)}
\enqa
Now as $|v(j)|\leq |v(k)|$,
\beqa
\left|\frac{\sum_{j\neq k}a_{kj}v(j)}{v(k)}\right|\leq\sum_{j\neq k}|a_{kj}|\left|\frac{v(j)}{v(k)}\right|\leq \sum_{j\neq k}|a_{kj}|=R_k.
\enqa
In other words $|\lambda-a_{kk}|\leq R_k$ $\bullet$
\end{proof}

\bigskip

Note that the diagonal entries of a stochastic operator driven by $\nu\in M_p(G)$ are all $\nu(e)$. Hence the radii, $R_t$, are all equal to $1-\nu(e)$. The theorem says, for any eigenvalue of the stochastic operator, $\lambda$, $|\lambda-\nu(e)|\leq 1-\nu(e)$.  Note $\tr P=|G|\nu(e)$. If $P$ is put in Jordan form, since the trace is basis independent, it is found that $\tr P=\sum_t \lambda_t$. Hence the average of the  eigenvalues\footnote{it would be interesting to try to apply this to obtain a bound for $\|\nu^{\star k}-\pi\|$} is equal to $\nu(e)$. Therefore, as $1$ is an eigenvalue, there are eigenvalues less than $\nu(e)$, i.e. $\lambda\leq \nu(e)$.  In the symmetric case, therefore, $\lambda_{|G|}- \nu(e)\leq0$ so that $-\lambda_{|G|}+\nu(e)\leq 1-\nu(e)$, thus $\lambda_{|G|}\geq-1+2\nu(e)$. In the general, not-necessarily-symmetric case, the eigenvalues are not necessarily real. However, with $|\lambda-\nu(e)|\leq 1-\nu(e)$, if $\nu(e)>1/2$, then the eigenvalues are bounded away from zero so that the stochastic operator, $P$, is invertible.
\newpage
\section{Comparison Techniques}
Whilst some random walks yield easily to analysis, others do not. There are a number of techniques, due to Diaconis \& Saloff-Coste \cite{chensixteen}, however, that allow comparison with a  simpler walk. Often the continuous analogue of a discrete random walk yields readily to analysis.  Diaconis \& Saloff-Coste \cite{chensixteen} present, in the symmetric case, the most general relationship between the discrete and  continuous time version of a given random walk. This paper also uses Dirichlet forms and the Courant minimax principle to estimate eigenvalues on a complicated walk from a simpler version.

\section{Lower Bounds}
The definition of variation distance immediately gives a technique for generating lower bounds. Given a test set $B\subset G$, immediately $\|\nu^{\star k}-\pi\|\geq |\nu^{\star k}(B)-\pi(B)|$. A very simple application uses the fact that $|\text{supp}(\nu^{\star k})|\leq |\Sigma|^k$. Let $A_k\subset G$ be the set where $\nu^{\star k}$ vanishes. Clearly
 \beqa
 |\stk(A_k)-\pi(A_k)|=\pi(A_k)\geq\frac{1}{|G|}(|G|-|\Sigma|^k)=1-\frac{|\Sigma|^k}{|G|}.
 \enqa
 Another elementary method for generating a lower bound using a test function is apparent via \beq\|\mu-\nu\|=\half\max_{\|\phi\|\leq 1}\left|\sum_{t\in G}(\mu(t)-\nu(t))\phi(t)\right|\label{lowbnd}.\enq
The discussion in Section \ref{charge} implies that if the (right) eigenvector $v_s$ is normalised to have $\|v_s\|_\infty=1$, then $v_s$ will have expectation zero under the random distribution, $\sum_t \pi(t)v_s(t)=0$.
 \wxn
\subsection{Proposition}
\wxo
 \emph{Let $u_\lambda$ be a real left eigenvector with eigenvalue $\lambda\neq 1$ and normalised such that $\pi+ u_\lambda\in M_p(G)$. Then
$$\|\nu^{\star k}-\pi\|\geq \half \|u_\lambda\|_1|\lambda|^k$$}
\begin{proof}
Let $\theta=\pi+u_\lambda$. Using the fact that a non-Dirac initial distribution $\theta$ converges faster than any Dirac measure (see (\ref{slow})), it is clear that $\|\theta P^k-\pi\|$ is a lower bound for $\|\nu^{\star k}-\pi\|$;
\beqa
\|\theta P^k-\pi\|=\|\pi+\lambda^ku_\lambda-\pi\|=\half \|u_\lambda\|_1 |\lambda|^k\,\,\,\bullet
\enqa
\end{proof}

\section{Volume \& Diameter Bounds}
By elegant analysis of the properties of the geometry of the random walk, bounds may be put on the  eigenvalues of $P$ and applications of the bounds of this Chapter give bounds on the variation distance. The geometry of the random walk is determined by its \textit{Cayley graph}.
   Suppose that $\xi$ is a random walk on $G$  with driving probability supported on  a generating set $\Sigma$. The Cayley graph  of the random walk is a directed graph with vertex set   identified with $G$.
 For any $g\in G$, $\sigma\in\Sigma$, the vertices corresponding to the elements $g$ and $\sigma g$ are joined by a directed edge. Thus the edge set  consists of pairs of the form $(g,\sigma g)$. The \textit{growth function} of the random walk is  $V(k):=|\Sigma^k|$ and the \textit{diameter} of $\xi$, $\Delta$,  is the minimum $k$ such $V(k)=|G|$. Say a random walk has \textit{$(A,d)$ moderate growth} if
\beq
\frac{V(k)}{V(\Delta)}\geq\frac{1}{A}\left(\frac{k}{\Delta}\right)^d\,,\,\,\,1\leq k\leq\Delta.
\enq
The following theorem appears in Diaconis \& Saloff-Coste \cite{chenseventeen}. The proof --- via the heavy machinery of path analysis, flows, two particular quadratic forms and some functional analysis --- is omitted. More details are to be found in \cite{chensixteen}. First an attractive lemma:
\wxn\subsection{Lemma}\wxo
\emph{Let $\xi$ be a symmetric random walk with diameter $\Delta$. Let $L:=\min\{\nu(s):s\in\Sigma\}$. Then, where $\lambda_2$ is the second largest eigenvalue\footnote{i.e. not necessarily $\lambda_\star$}}
\beq
\lambda_2\leq 1-\frac{L}{\Delta^2}
\enq
\newpage
\wxn\subsection{Theorem\label{12}}\wxo
\emph{Let $\xi$ be a symmetric random walk with $(A,d)$ moderate growth. Then for $k=(1+c)\Delta^2/L$, with $c>0$:
\beq
\|\nu^{\star k}-\pi\|\leq B e^{-c}
\enq
where $B=2^{d(d+3)/4}\sqrt{A}$.
\new Conversely, for $k=c\Delta^2/(2^{4d+2}A^2)$:
\beq
\|\nu^{\star k}-\pi\|\geq \half e^{-c}
\enq}
\subsubsection{Example: The Heisenberg Group}
Consider the set of matrices:
\beq
H_3(n)=\left(\begin{array}{ccc}
1 & a & b
\\0 & 1 & c
\\ 0 & 0 & 1
\end{array}\right)
\enq
where $a$, $b$, $c\in\Z_n$. With matrix multiplication modulo $n$, $H_3(n)$ forms a group of order $n^3$. The random walk driven by the measure $\nu_n\in M_p(H_3(n))$ constant on the matrices $(a,b,c)=(\pm1,0,0)$, $(0,0,\pm1)$, $(0,0,0)$ is ergodic. Diaconis \& Saloff-Coste \cite{chenseventeen} have shown that the random walk has diameter $n-1\leq \Delta\leq n+1$ and volume growth function $V(k)\geq k^3/6$ ($1\leq k\leq n+1$). With order\footnote{$n^3\leq 8(n-1)^3\leq 8\Delta^3$ for $n\geq 4$.} $|H_3(n)|\leq 8\Delta^3$,
$$\frac{V(k)}{V(\Delta)}\geq \frac{k^3/6}{8\Delta^3}=\frac{1}{48}\left(\frac{k}{\Delta}\right)^3\text{ , for }1\leq k\leq\Delta,$$  the random walk has $(48,3)$ moderate growth.
 \new Precise application of Theorem \ref{12} yields for constants $A$, $A^\prime$, $B$, $B^\prime$:
\beq
A^\prime e^{-B^\prime k/n^2}\leq\|\nu_n^{\star k}-\pi\|\leq A e^{-Bk/n^2}
\enq
Hence order $n^2$ steps are necessary for convergence to random.

%

\chapter{Diaconis-Fourier Theory}
Much of the precluding analysis passes neatly into the case of the classical Markov theory for a random walk on a finite set $X$. It has been seen that this analysis culminates in the result that the rate
of convergence of to a stationary state is related heavily to the
second largest eigenvalue of the stochastic operator. As a rule the
calculation of the second highest eigenvalue is too cumbersome for
larger groups and further the bound is not particularly sharp due to
the information loss in disregarding the rest of the spectrum of the stochastic operator.

 In his seminal monograph \cite{PD}, Diaconis utilises the group structure to produce bounds for rates of convergence. He uses Fourier methods and representation
theory to produce bounds  that are invariably sharper as the entire spectrum is utilised.  This chapter follows his approach.

\section{Basics of Representations and Characters}
A \textit{representation} $\rho$
of a finite group $G$ is a group homomorphism from $G$ into  $GL(V)$ for some vector space $V$. The
dimension of the vector space\footnote{at this point the underlying vector space may be infinite dimensional but later it will be seen that the only representations of any interest are of finite dimension. Also the underlying field is unspecified at this point but  later it will be seen that the only representations of any interest will be over complex vector spaces.} is called the \textit{dimension of
$\rho$} and is denoted by $d_\rho$.   If $W$ is a subspace of $V$
invariant under $\rho(G)$, then
$\rho_{|W}$ is called a \textit{subrepresentation}.
 \newpage If $\la.,.\ra$ is an inner product on $V$,  $\la u,v\ra_\rho=\sum_t\la\rho(t)u,\rho(t)v\ra$ defines another, and further the orthogonal complement of $W$ with respect to $\la.,.\ra_\rho$, $W^\perp$, is also invariant under $\rho$. Hence, every representation splits into a direct sum of subrepresentations. Both
$\{\bld{0}\}$ and $V$ itself yield trivial subrepresentations. A
representation $\rho$ that admits no non-trivial subrepresentations
is called \textit{irreducible}. An example of an irreducible representation is the \textit{trivial representation}, $\tau$, which maps $G$ to 1: $\rho(s)z=z$, $z\in\C$. Inductively, therefore, every representation is a direct sum of irreducible
representations. A quick calculation shows $\la \rho(s)u,\rho(s)v\ra_\rho=\la u,v\ra_\rho$, hence $\|u\|_\rho=\|\rho(s)v\|_\rho$ so the operators  $\rho(s)$ are isometries and are thus unitary. Two representations, $\rho$ acting on $V$ and $\varrho$ acting on $W$;
 are \textit{equivalent as representations}, $\rho\equiv \varrho$, if there is a bijective
linear map $f\in L(V,W)$  such that $\varrho\circ f=f\circ
\rho$. In this context $f$ is said to \textit{intertwine} $\varrho$ and $\rho$.
%

\subsubsection{Example: A Two Dimensional Representation of the Dihedral Group}
The dihedral group $D_4$, the group of symmetries of the square, admits a natural representation $\rho$. The elements of $D_4$ are the rotations $r_{0}$, $r_{\pi/2}$, $r_{\pi}$, $r_{3\pi/2}$ and reflections $(12)$, $(13)$, $(14)$, $(23)$. If the vertices of the square are inscribed in a unit circle at the poles\footnote{i.e. the coordinates $(\pm1,0)$, $(0,\pm1)$.} then $\rho(r_\theta)$ are the rotation matrices:
\beqa
\rho(r_\theta)=\left(\begin{array}{cc}\cos\theta & -\sin\theta\\ \sin\theta & \cos \theta\end{array}\right)
\enqa
Similarly the reflections have action as reflection in $y=x$, $y=-x$, $y=0$ and $x=0$ which have matrix representations:
\beqa
\begin{array}{cc}
\rho((12))=\left(\begin{array}{cc}1 & 0\\ 0 & -1\end{array}\right) & \rho((13))=\left(\begin{array}{cc}-1 & 0\\ 0 & 1\end{array}\right)
\\ \rho((14))=\left(\begin{array}{cc}0 & 1\\ 1 & 0\end{array}\right) & \rho((23))=\left(\begin{array}{cc}0 & -1\\ -1 & 0\end{array}\right)
\end{array}
\enqa

\newpage

\wxn\subsection{Schur's Lemma}\wxo
\emph{Let $\rho_1:G\raw GL(V_1)$ and $\rho_2:G\raw GL(V_2)$ be two irreducible representations of $G$, and let $f\in L(V_1,V_2)$  be an intertwiner. Then
\begin{enumerate}
\item If $\rho_1$ and $\rho_2$ are not equivalent $f\equiv0$.
\item If $V_1=:V:=V_2$ is complex, and $\rho_1:=\rho=:\rho_2$, $f=\lambda I$, for some, $\lambda\in\C$.
\end{enumerate}}
\begin{proof}
The straightforward calculations $f(\rho_1(G)\ker f)=\rho_2(G)f(\ker f)=0$ and\new $\rho_2(G)\im f=f(\rho_1(G)V_1)$ show that $\ker f$ and $\im f$ are invariant subspaces. By irreducibility both the kernel and image of $f$ are trivial or the whole space.
\begin{enumerate}
\item Suppose $f\not\equiv 0$. Hence $\ker f=\{0\}$ and $\im f=V_2$ so $f$ is an isomorphism as it is linear. However this would imply that $\rho_1$ and $\rho_2$ are equivalent as representations, a contradiction. Thence $f\equiv 0$.
\item If $f\equiv 0$ then $f=0.I$.  Suppose again $f\not\equiv 0$. Then $f$ has a non-zero eigenvalue $\lambda\in \C$ with associated non-zero eigenvector $v_\lambda\neq0$. Let $f_\lambda=f-\lambda I$. A quick calculation shows that $\rho(G)f_\lambda(V)=f_\lambda(\rho(G)V)$, hence $f_\lambda$ is an intertwiner. Note that $\ker f_\lambda\neq\{0\}$ as $v_\lambda\in\ker f_\lambda$. Thence $\ker f_\lambda =V$, that is $f_\lambda\equiv 0$, which implies $f=\lambda I\,\,\,\bullet$
\end{enumerate}
\end{proof}

\bigskip

Let $\rho_1:G\raw GL(V_1)$ and $\rho_2:G\raw GL(V_2)$ be two irreducible representations of $G$ and $h_0\in L(V_1,V_2)$. Let \beq h=\frac{1}{|G|}\sum_{t\in G}\rho^{-1}_2(t)h_0\rho_1(t)\enq
A quick verification shows that $h$ is an intertwiner of $\rho_1$ and $\rho_2$, and by recourse to Schur's Lemma $h\equiv 0$ in the case where $\rho_1\not\equiv\rho_2$, and $h=\lambda I$ when $\rho_1\equiv\rho_2$. In the case $\rho_1\equiv \rho_2$, taking traces gives $\lambda=\tr h/d_\rho$ and a further calculation shows $\tr h=\tr h_0$.
Suppose $\rho_1$ and $\rho_2$ are given in matrix form as $\rho_1(s)=(r^1_{ij}(s))$ and $\rho_2(s)=(r^2_{ij}(s))$. The linear maps $h$ and $h_0$ are defined by matrices ${x}_{ij}$ and ${x}^0_{ij}$. In particular,
\beq
x_{ij}=\frac{1}{|G|}\sum_{\underset{\lambda,\mu}{t\in G}} r^2_{i\lambda}(t^{-1}) {x}^0_{\lambda \mu}r^1_{\mu j}(t)\label{hh0}
\enq
Suppose $\rho_1\not\equiv\rho_2$ so that $h\equiv0$ when defined by $h_0=\delta_{kl}$. In this case $x_{ij}=0$ and (\ref{hh0}) collapses to
\beq
\frac{1}{|G|}\sum_{t\in G}r^2_{ik}(t^{-1})r^1_{lj}(t)=0\,\,,\,\forall\, i,k,l,j.\label{o1}
\enq
In the case where $\rho_1\equiv\rho_2$, $h=\lambda I$, where, in matrix elements, $\lambda=\sum_{m}x^0_{mm}/d_\rho$. When $h$ is defined by $h_0=\delta_{kl}$, (\ref{hh0}) collapses to
\beq
\frac{1}{|G|}\sum_{t\in G}r^2_{ik}(t^{-1})r^1_{lj}(t)=\frac{\delta_{ij}\delta_{kl}}{d_\rho}.\label{o2}
\enq
Note again that $\rho(s)$ is a unitary operator so that $\rho(s)^\star=\rho^{-1}(s)$, thence $\overline{r_{ji}(s)}=r_{ij}(s^{-1})$. A simple rearrangement of (\ref{o1}) and (\ref{o2}) using this fact show that the matrix elements of the irreducible representations are orthogonal in $F(G)$.

\bigskip

If $\rho$ is a representation,  the \textit{character} of
$\rho$, $\chi_\rho(s):=\tr \rho(s)$. Using the preceding remarks, it can be shown that the characters of the irreducible representations are orthonormal in $F(G)$. If $\rho_1$ and $\rho_2$ are representations with characters $\chi_1$ and $\chi_2$, by choosing a basis so that the matrix of $\rho_1\oplus \rho_2$ is a
block $2\times2$ matrix with $\rho_1$ in the $(1,1)$ position and
$\rho_2$ in the $(2,2)$ position, taking traces shows that the character of $\rho_1\oplus\rho_2$ is $\chi_1+\chi_2$. Suppose now $\rho$ is a representation with character $\phi$ that decomposes into a direct sum of irreducible representations $\rho=\rho_1\oplus\cdots\oplus \rho_k$. If each of the $\rho_i$ have character $\chi_i$, then $\phi=\chi_1+\cdots+\chi_k$. If $\rho^\prime$ is an irreducible representation with character $\chi$, then $(\phi|\chi)=\sum_i(\chi_i|\chi)$. By orthonormality, $(\chi_i|\chi)=$ 0 or 1 as $\chi_i$ is, or is not,
equivalent to $\chi$. Thence, the number of $\rho_i$ equivalent to $\rho^\prime$ equals $(\phi|\chi)$.

\bigskip

A canonical representation is the \textit{regular representation}; defined with respect to a complex
vector space with basis $\{e_s\}$ indexed by $s\in G$ via $r(s)(e_t):=e_{st}$.
Observe that the underlying vector space is isomorphic to $F(G)$. It is a simple exercise to show that $\chi_r(e)=|G|$ and zero elsewhere. This implies that for an irreducible representation $\rho_i$, $(\chi_r|\chi_i)=\chi_i(e)^\star=\overline{\tr I_{d_i}}=d_i$ so that $\chi_r(s)=\sum_i d_i\chi_i(s)$, where the sum is over all irreducible representations. Letting $s=e$ here yields $\sum_i d_i^2=|G|$. Now it can be seen that the matrix entries of the irreducible representations form an orthogonal basis for $F(G)$ because they are orthogonal and there are $\sum_i d_i^2=|G|$ of them: $\text{dim}\,F(G)=|G|$.

\section{Fourier Theory}
Let $f\in F(G)$ and $\rho$ a representation of $G$.
The \textit{Fourier Transform of $f$ at the representation $\rho$}
is the operator $\widehat{f}(\rho)=\sum_s f(s)\rho(s)$. This Fourier transform satisfies an inversion theorem, a Plancherel Formula; and, of course, a Convolution Theorem $\widehat{f\star h}(\rho)=\widehat{f}(\rho)\widehat{h}(\rho)$ whose proof is rudimentary.
\wxn\subsection{Fourier Inversion Theorem\label{fit}}\wxo
\emph{Let $f\in F(G)$, then, where the sum is over
irreducible representations,} \beq f(s)=\frac{1}{|G|}\sum_i d_i\,\tr
(\rho_i(s^{-1})\widehat{f}(\rho_i)) \enq
\begin{proof}
Both sides are linear in $f$ so it is sufficient to check the
formula for $f= \delta_t$. Then
$\widehat{f}(\rho_i)=\rho_i(t)$, and the right hand side equals
$$\frac{1}{|G|}\sum_id_i\tr (\rho_i(s^{-1})\rho_i(t))=\frac{1}{|G|}\sum_i d_i\chi_i(s^{-1}t)=\frac{1}{|G|}\chi_r(s^{-1}t)$$
When $s=t$ this equals 1; otherwise it is 0; i.e. it equals $\delta_t$ $\bullet$
\end{proof}

\wxn\subsection{Plancherel Formula\label{pf}}\wxo
\emph{Let $f$, $h\in F(G)$, then}
\beq
\sum_{s\in G} f(s^{-1})h(s)=\frac{1}{|G|}\sum_i
d_i\,\tr(\widehat{f}(\rho_i)\widehat{h}(\rho_i))
\enq
\begin{proof}
Both sides are linear in $f$; so consider $f= \delta_t$. Using the Fourier Inversion Theorem
\beqa
h(t^{-1})=\sum_{s\in G}\delta_{t}(s^{-1})h(s)=\frac{1}{|G|}\sum_i d_i\,\tr(\rho_i(t)\widehat{h}(\rho_i))
\enqa

However, $\rho_i(t)$ is nothing but $\widehat{\delta}_{t}(\rho_i)$ so the formula is verified $\bullet$
\end{proof}

\newpage

In the sequel, mostly elements of $M_p(G)$ viewed as elements of $F(G)$ are considered. Let $\mu\in M_p(G)$ and let $\check{\mu}(s):=\mu(s^{-1})$. After a reindex, $t\mapsto t^{-1}$, $\widehat{\check{\mu}}(\rho)=\sum_t\mu(t)\rho(t^{-1})$. With the unitary nature of the representation, and the fact that $\mu=\overline{\mu}$ as $\mu\in \R$, in fact $\widehat{\check{\mu}}(\rho)=\widehat{\mu}(\rho)^\star$.  Hence, for $\nu\in M_p(G)$:
\beq
 \sum_{t\in G}\mu(t)\nu(t)=\frac{1}{|G|}\sum d_i\,\tr\left[\widehat{\nu}(\rho_i)\widehat{\mu}(\rho_i)^\star\right]\label{37}
\enq
With the aid of two quick facts the celebrated Upper Bound Lemma of Diaconis and Shahshahani \cite{PD,dash} may be proven. The first of these is the straightforward calculation that for all $\nu\in M_p(G)$, at the trivial representation $\tau$, $\widehat{\nu}(\tau)=\sum_t\nu(t)=1$. The second comprises a lemma.
\wxn\subsection{Lemma\label{l2}}\wxo
\emph{At a non-trivial irreducible representation, $\rho$, the Fourier transform of the random distribution, $\pi$, vanishes: $\widehat{\pi}(\rho)=\bld{0}$.}
\begin{proof}
First note that $h=\sum_{t\in G}\rho(t)$ is a linear map, invariant under any $\rho(s)$: $\rho(s)h=h=h\rho(s)$. As a consequence both $\ker h$ and $\im h$ are invariant subspaces. By irreducibility, both the kernel and the image of $h$ are trivial or the whole space. Suppose $\ker h=\{\bld{0}\}$ and $\im h=V$. For any $v\in V$, $\rho(s)hv=hv$. Hit both sides with $h^{-1}$: $h^{-1}\rho(s)hv=v$. Now use the fact that $\rho(s)$ and $h$ commute to show $\rho(s)v=v$. Hence $\rho$ is trivial. Therefore $\ker h=V$, $\im h=\{\bld{0}\}$, i.e. $h=\bld{0}$. Now $\widehat{\pi}(\rho)=\sum_t\pi(t)\rho(t)=h/|G|=\bld{0}\,\,\,\bullet$
\end{proof}

\wxn\subsection{Upper Bound Lemma}\wxo
\emph{Let $\nu$ be a probability on a finite group $G$. Then} \beq
\|\nu-\pi\|^2\leq\frac{1}{4}\sum_i
d_i\,\tr(\widehat{\nu}(\rho_i)\widehat{\nu}(\rho_i)^\star), \enq \emph{where the
sum is over all non-trivial irreducible representations}.
\begin{proof}
Using the Cauchy-Schwarz Inequality
\begin{align*} 4\|\nu-\pi \|^2&=\left\{\sum_{t\in G}|\nu(t)-\pi (t)|\right\}^2
\\&\leq |G|\sum_{t\in G}|\nu(t)-\pi (t)|^2= |G|\sum_{t\in G}(\nu-\pi) (t)(\nu-\pi) (t),
\end{align*}
where of course $\nu-\pi$ is a real function.  Thus, by (\ref{37})
\beqa
 4\|\nu-\pi \|^2\leq \cancel{|G|}\frac{1}{\cancel{|G|}}\sum_i d_i\, \tr\left[\widehat{(\nu-\pi )}(\rho_i)\widehat{(\nu-\pi )}(\rho_i)^\star\right]
\enqa
Now $\widehat{(\nu-\pi )}(\rho)=\sum_{t}(\nu(t)-\pi (t))\rho(t)=\widehat{\nu}(\rho)-\widehat{\pi }(\rho)$. With the preceding facts:
\beqa
\widehat{(\nu-\pi )}(\rho)=\left\{\begin{array}{cc} \bld{0} & \text{ if $\rho$ is trivial}\\\widehat{\nu}(\rho) & \text{ if $\rho$ is non-trivial and irreducible}\end{array}\right.
\enqa
So  therefore
\beqa 4\|\nu-\pi \|^2\leq \sum_i d_i \,\tr(\widehat{\nu}(\rho_i)\widehat{\nu}(\rho_i)^*),
\enqa
where the sum is over all non-trivial representations$\,\,\,\bullet$
\end{proof}

\bigskip

%

This bounds are applicable to $\|\nu^{\star k}-u\|$ via the
Convolution Theorem: $\widehat{\nu^{\star k}}(\rho)=\widehat{\nu}(\rho)^k$.

\section{Number of Irreducible Representations}
Let $G$ be a group and $g,h$ elements of $G$. An element $g\in G$ is \textit{conjugate} to $h$, $g\sim h$, if  there exists $ t\in G$ such that $h=tgt^{-1}$. Conjugacy is an equivalence relation on a group \cite{H}, and hence forms a partition of $G$ into disjoint \textit{conjugacy classes} $G=[s_1]_\sim\cup[s_2]_\sim\cup\cdots\cup[s_r]_\sim$, where
\beq
[s]_\sim=\{g\in G:\exists\,t\in G,g=tst^{-1} \}= \{tst^{-1}:t\in G\}.
\enq
A complex function $f\in F(G)$ is a \textit{class function} if for all conjugacy classes $[s_i]_{\sim}\subset G$, $f_{|{[s_i]}_\sim}=\lambda$, for some $\lambda\in \C$. Let $\mathfrak{Cl}(G)$ be the subspace of $F(G)$ consisting of all class functions. The characters of a representation are class functions. Let $f\in\mathfrak{Cl}(G)$ and  $\rho$ be an irreducible representation. Note that $\rho(s)\widehat{f}(\rho)\rho(s^{-1})=\sum_tf(t)\rho(sts^{-1})$, and with a reindexing $t\mapsto s^{-1}ts$, it is clear that $f$ is an intertwiner for $\rho$. Thus, by Schur's Lemma, $\widehat{f}(\rho)=\lambda I$. Taking traces gives $\lambda=\tr(\widehat{f}(\rho))/d_\rho=\sum_tf(t)\chi(t)/d_\rho=|G|(f|\chi^\star)/d_\rho$.

\wxn\subsection{Theorem\label{th6}}\wxo
\emph{The characters of the irreducible representations $\chi_1,\chi_2,\dots,\chi_l$ form an orthonormal basis for $\mathfrak{Cl}(G)$.}
\begin{proof}
Characters are orthonormal class functions. As $\mathfrak{Cl}(G)$ together with $(.|.)$ forms an inner product space, and $\Omega=\text{span}\{\chi_i\}$ is a subspace: $\mathfrak{Cl}=\Omega\oplus\Omega^\perp$. Let $f\in \mathfrak{Cl}$ have the decomposition $f=g+h$, with $g\in\Omega$, $h\in \Omega^\perp$.  Therefore for all irreducible representations $\chi_i$: $(h|\chi_i^\star)=0$. The preceding remarks indicate that $\widehat{h}(\rho)=|G|(h|\chi_\rho^\star)I/d_\rho=\bld{0}$. The Fourier Inversion Theorem yields:
\beqa
h(s)=\frac{1}{|G|}\sum_id_i\tr\left(\rho_i(s^{-1})\widehat{h}(\rho_i)\right)\equiv0.
\enqa
Hence therefore $\Omega^\perp=\{\bld{0}\}$ and the characters of the irreducible representations span $\mathfrak{Cl}(G)$.
$\,\,\,\bullet$
\end{proof}

\wxn\subsection{Theorem}\wxo
\emph{The number of irreducible representations equals the number of conjugacy classes.}
\begin{proof}
Theorem \ref{th6} gives the number of irreducible representations, $l$:
\beqa
l=\text{dim}(\mathfrak{Cl}(G))
\enqa
A class function can be defined to have an arbitrary value on each conjugacy class, so $\text{dim}(\mathfrak{Cl}(G))$ is the number of conjugacy classes$\,\,\,\bullet$
\end{proof}

\bigskip

As an immediate corollary, all the irreducible representations of an Abelian group $G$ have degree 1. To see this note if $G$ is Abelian, there are $|G|$ conjugacy classes, so $|G|$ terms in the sum $\sum_id_i^2=|G|$, each of which must be 1. Hence if $G$ has $l$ conjugacy classes and $l$ representations are found, if the $l$ representations are inequivalent and irreducible, all the irreducible representations have been found.

\wxn\subsection{Theorem}\wxo
\emph{Two irreducible representations with the same character are equivalent.}
\begin{proof}
 Suppose $\chi_1$, $\chi_2$ are identical characters of non-equivalent irreducible representations $\rho_1$ and $\rho_2$,
 \beqa
 (\chi_1|\chi_2)=\frac{1}{|G|}\sum_{t\in G}\chi_1(t)\overline{\chi_1(t)}\underset{\chi_1(e)\neq0}{>}0.
 \enqa
 However the characters of irreducible representations are orthonormal. This is a contradiction; hence $\rho_1\equiv\rho_2\,\,\,\bullet$
\end{proof}

 \wxn\subsection{Theorem\label{cha}}\wxo
 \emph{Let $\chi$ be the character of a representation $\rho$, then $\rho$ is an irreducible representation if and only if $(\chi|\chi)=1$.}
\begin{proof}
Clearly if $\rho$ is irreducible $(\chi|\chi)=1$. Suppose for the converse that $(\chi|\chi)=1$. Any representation $\rho$ is the direct sum of irreducible representations $\{\rho_i\}$ with character $\chi=\chi_1+\chi_2+\cdots+\chi_m$. Therefore if $(\chi|\chi)$ must equal $1$, then there exists a unique $\rho_k$ such that $\rho\equiv\rho_k$$\,\,\,\bullet$
\end{proof}

\subsubsection{Example: The Quaternion Group, $\mathcal{Q}$}
Consider the quaternion group $\mathcal{Q}=\{\pm1,\pm i,\pm j,\pm k\}$ where $1$ is the identity. Multiplication in $\mathcal{Q}$ is defined by $(-1)^2=1$ and $i^2=j^2=k^2=ijk=-1$, where $-1$ commutes with everything. The quaternion group has five conjugacy classes $\{1\}$, $\{-1\}$, $\{\pm i\}$, $\{\pm j\}$ and $\{\pm k\}$ and thus five irreducible representations.  As $\sum_i d_i^2=|G|$, the there must be one irreducible representation of degree 2 and four of degree 1. Consider the linear map $\rho:\mathcal{Q}\raw GL(\C^2)$ given by:
\beq
\begin{array}{cc}
\rho(i)=\left(\begin{array}{cc}i & 0\\ 0 & -i\end{array}\right) & \rho(k)=\left(\begin{array}{cc}0 & -1\\ 1 & 0\end{array}\right)
\\ \rho(j)=\left(\begin{array}{cc}0 & i\\ i & 0\end{array}\right) & \rho(1)=I,\,\,\rho(-s)=-\rho(s)
\end{array}
\enq

Straightforward calculations show that $\rho$ is a representation. Also  $(\chi|\chi)=1$, and in light of Theorem \ref{cha}, $\rho$ is the two dimensional irreducible representation. Let $\tau:\mathcal{Q}\raw GL(\C)$ be the trivial representation; it is the second irreducible representation. Let $\rho_i:\mathcal{Q}\raw GL(\C)$ (respectively $\rho_j$, $\rho_k$) be defined by:
\beq
\rho_i(s):=\left\{\begin{array}{cc}1  & \text{ if }s\in\la i \ra\\ -1 & \text{ if }s\not\in\la i \ra\end{array}\right.
\enq
This is a one-dimensional representation  so is irreducible. It is an easy calculation to show that $\{\tau,\chi_i,\chi_j,\chi_k\}$ is an orthogonal set so comprise four inequivalent representations. Hence the set of irreducible representations of $\mathcal{Q}$ are given by $\{\rho,\tau,\rho_i,\rho_j,\rho_k\}$.
\section{Simple Walk on the Circle}
Consider the walk on $\{\Z_n,\oplus\}$  driven by
\beq
\nu_n(s):=\left\{\begin{array}{cc}
\half & \text{ if }s=\pm1
\\[1ex] 0 & \text{ otherwise}
\end{array}\right.
\enq

$\Z_n$ is an Abelian group, so all irreducible representations have degree 1. Any $\rho$ is determined by the image of $1$: $\rho(s)=\rho(1^s)=\rho(1)^s$. Also $1^n=0$, hence $\rho(1)^n=\rho(1^n)=\rho(0)=1$
so $\rho(1)$ must be a $n$-th root of unity. There are $n$ such: $e^{2\pi i t/n}$, $t=0,1,2,\dots,n-1$. Each gives a representation $\rho_t(s)=e^{2\pi i t s/n}$.  Now some results used in the Lower Bound; see Appendix A for proof.

\wxn\subsection{Lemma\label{zpl}}\wxo
\emph{The following (in)equalities hold.}
\begin{enumerate}
\item \emph{For any odd $n$ and $k\in\N$},
\beq
\sum_{t=1}^{n-1} \cos^{2k}(2\pi t/n)=2\sum_{t=1}^{(n-1)/2}\cos^{2k}(\pi t/n)\label{lem1}
\enq
\item \emph{For $x\in[0,\pi/2]$},
\beq
\cos x\leq e^{-x^2/2}\label{lem2}
\enq
\item \emph{For any $x>0$}
\beq
\sum_{j=1}^\infty e^{-(j^2-1)x}\leq \sum_{j=0}^\infty e^{-3jx}\label{lem3}
\enq
\item \emph{For $x\in[0,\pi/6]$},
\beq
\cos x\geq e^{-x^2/2-x^4/2}\label{lem4}
\enq
\end{enumerate}
\wxn\subsection{Upper and Lower Bounds\label{circ}}\wxo
\emph{For $k\geq n^2/40$, with $n$ odd,
\beq
\|\nu_n^{\star k}-\pi_n \|\leq e^{-\pi^2 k/2n^2}
\enq
Conversely, for $n\geq 7$, and \text{any} $k$
\beq
\|\nu_n^{\star k}-\pi_n\|\geq\half e^{-\pi^2 k/2n^2-\pi^4 k/2n^4}.
\enq}
\begin{proof}
The Fourier transform of $\nu_n$ at $\rho_s$ is:
\beqa
\widehat{\nu_n}(\rho_s)=\sum_{t=0}^{n-1}\nu_n(t)e^{2\pi i s t/n}=\frac{1}{2}\,e^{2\pi i s/n}+\half\,e^{-2\pi i s/n}=\cos\left(\frac{2\pi s}{n}\right).
\enqa
The Upper Bound Lemma and (\ref{lem1}) yield
\beqa
\|\nu_n^{\star k}-\pi_n \|^2\leq \frac{1}{4}\sum_{t=1}^{n-1}\cos^{2k}\left(\frac{2\pi t}{n}\right)=\half \sum_{t=1}^{(n-1)/2}\cos^{2k}\left(\frac{\pi t}{n}\right).
\enqa
Applying (\ref{lem2}) yields
\beqa
\|\nu_n^{\star k}-\pi_n \|^2\leq\half \sum_{t=1}^{(n-1)/2}e^{-\pi^2t^2k/n^2}\leq \half\,e^{-\pi^2k/n^2}\sum_{t=1}^\infty e^{-\pi^2 (t^2-1)k/n^2},
\enqa
and so with (\ref{lem3})
\beqa
 \|\nu_n^{\star k}-\pi_n \|^2\leq\half\, e^{-\pi^2k/n}\sum_{t=0}^\infty e^{-3\pi^2 t k/n^2}= \half\,\frac{e^{-\pi^2 k/n^2}}{1-e^{-3\pi^2 k/n^2}}.
\enqa
Now since $k\geq n^2/40$, $2\left(1-e^{-3\pi^2k/n^2}\right)>1$, and
it follows that
\beqa
\|\nu_n^{\star k}-\pi_n \|\leq e^{-\pi^2 k/2n^2}
\enqa

\bigskip

For the lower bound, consider the norm 1 function $\phi(s)=\rho_{\bar{s}}(s)=\cos(2\pi s \bar{s}/n)$ where $\bar{s}=(n-1)/2$. By Lemma \ref{l2}, $\phi(s)$ has zero expectation under the random distribution.
Now an application of (\ref{lowbnd}) gives
\beqa
\|\nu_n^{\star k}-\pi_n\|\geq \half \left|\sum_{t\in G}\nu_n^{\star k}(t)\phi(t)\right|=\half\left|\widehat{\nu_n^{\star k}}\right|=\half\left|\widehat{\nu_n}(\rho_{\bar{s}})\right|^k
\enqa
Now $\widehat{\nu_n}(\rho_{\bar{s}})=\cos (2\pi\bar{s}/n)=-\cos (\pi/n)$ by a quick calculation. By (\ref{lem4}), for $\pi/n\leq \pi/6$:
\beqa
\|\nu_n^{\star k}-\pi_n\|\geq \half\left|\cos\frac{\pi}{n}\right|^k\geq \half e^{-\pi^2 k/2n^2-\pi^4 k/2n^4}\,\,\,\bullet
\enqa
\end{proof}

\subsubsection*{Remark}
If $n$ is even then $\{1,-1\}$ lies in the coset of odd numbers of the normal subgroup $\{0,2,\dots,n-2\}=:H\lhd \Z_{n}$, and so the walk is not ergodic by Theorem \ref{eth}.

\begin{figure}[h]\cone\epsfig{figure=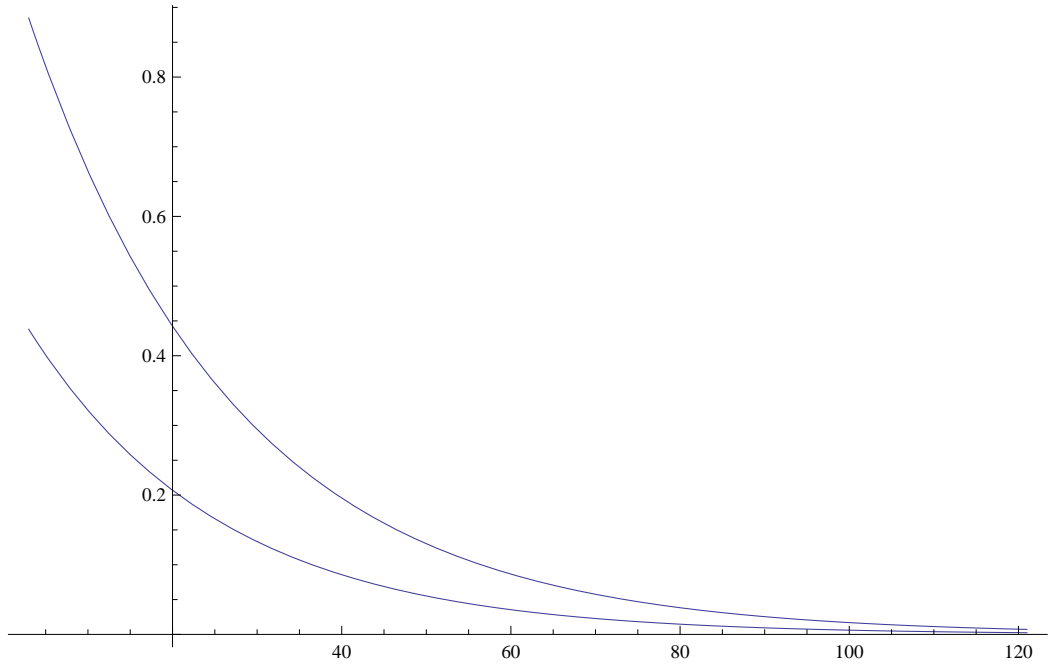}\caption{A plot of the upper and lower bound  for $n=11$.}\ctwo\end{figure}

\section{Nearest Neighbour Walk on the $n$-Cube}
Consider the walk on $\{\Z_2^n,\oplus^n_2\}$, $n>1$,  driven by
\beq
\nu_n(s):=\left\{\begin{array}{cc}
\frac{1}{n+1} & \text{ if }w(s)=0\text{ or }1
\\[2ex] 0 & \text{ otherwise}
\end{array}\right.
\enq
where $w(s)$, the weight of $s=(s_1,s_2,\dots,s_n)$, is given by the sum in $\N$:
\beq
w(s)=\sum_{i=1}^ns_i
\enq

$\Z_2^n$ is an Abelian group, so all irreducible representations have degree 1. It is a simple verification to show that each are given by $\rho_t(s)=(-1)^{t\cdot s}$.  Now some results used in the Upper Bound; see Appendix A for proof.
\wxn\subsection{Lemma\label{hcl}}\wxo
\emph{The following inequalities hold.}
\begin{enumerate}
\item \emph{If $l\leq n/2$},
\beq
{n \choose l}\left(1-\frac{2l}{n+1}\right)^{2k}\geq {n\choose n+1-l}\left(1-\frac{2(n+1-l)}{n+1}\right)^{2k}\label{lem5}
\enq
\item \emph{When $a\leq b$},
\beq
{a\choose b}\leq \frac{a^b}{b!}\label{lem6}
\enq
\item \emph{Let $n\in\N$, $c>0$. If $k=(n+1)(\log n+c)/4$}
\beq
\left(1-\frac{2j}{n+1}\right)^{2k}\leq e^{-j\log n-jc}\label{lem7}
\enq
\end{enumerate}
\wxn\subsection{Upper Bound}\wxo
\emph{For $k=(n+1)(\log n+c)/4$, $c>0$:}
\beq
\|\nu_n^{\star k}-\pi_n \|^2\leq \half\left(e^{e^{-c}}-1\right)
\enq
\begin{proof}
Let $\{e_i\}$ denote the standard basis\footnote{of the finite vector space $\Z_2^n$ with underlying field $\Z_2$.} of $\Z_2^n$:
\beqa
\widehat{\nu_n}(\rho_s)=\sum_{t\in \Z_2^n} (-1)^{s\cdot t}\nu_n(t)=\frac{1}{n+1}\left[1+\sum_{i=1}^n(-1)^{s\cdot e_i}\right]
\enqa
Now  $s\cdot e_i=s_i$ so
\begin{align*}
\widehat{\nu_n}(\rho_s)&=\frac{1}{n+1}\left[1+\sum_{i=1}^n(-1)^{s_i}\right]
\\ &=\frac{1}{n+1}\left[1+\sum_{s_i=1}(-1)+\sum_{s_i=0}(1)\right]
\\ &=\frac{1}{n+1}\left[1+w(s)(-1)+(n-w(s))(1)\right]
\\ &=\frac{n+1-2w(s)}{n+1}=1-\frac{2w(s)}{n+1}
\end{align*}
Thus Upper Bound Lemma gives (summing over weights on the right equality):
\beq
\|\nu_n^{\star k}-\pi_n \|^2\leq \frac{1}{4}\sum_{t\neq 0}{\widehat{\nu_n}(\rho_t)}^{2k}= \frac{1}{4}\sum_{j=1}^n {n\choose j}\left(1-\frac{2j}{n+1}\right)^{2k}.\label{ub1}
\enq
Let $n/2\leq j\leq n$ such that $j=n+1-l$ (i.e. $l\in\{1,2,\dots,\lfloor n/2\rfloor\}$) and consider the $(n+1-l)$th (i.e. $j$th) term in this sum. By (\ref{lem5}), the $l$th term dominates this term, and for $l\in\{1,2,\dots,\lfloor n/2\rfloor\}$,
\beq
{n\choose l}\left(1-\frac{2l}{n+1}\right)^{2k}+{n\choose n+1-l}\left(1-\frac{2(n+1-l)}{n+1}\right)^{2k}\leq 2 {n\choose l}\left(1-\frac{2l}{n+1}\right)^{2k}
\enq
Noting that the `middle' term (i.e. $n$ odd) is unaffected, (\ref{ub1}) is thus dominated by a sum of $\lceil n/2\rceil$ terms. Therefore, with (\ref{lem6})
\beqa
 \|\nu_n^{\star k}-\pi_n \|^2\leq \frac{1}{2}\sum_{j=1}^{\lceil n/2\rceil}\frac{n^j}{j!}\left(1-\frac{2j}{n+1}\right)^{2k}.
\enqa
Applying (\ref{lem7}) and noting $n^j=e^{j\log n}$,
\begin{align*}
\|\nu_n^{\star k}-\pi_n \|^2&\leq \half\sum_{j=1}^{\lceil n/2\rceil}\frac{e^{j\log n}}{j!}e^{-j\log n-jc}
 =\half\sum_{j=1}^{\lceil n/2\rceil}\frac{e^{-jc}}{j!}
\\ &\leq \half \sum_{j=1}^\infty \frac{(e^{-c})^j}{j!}=\half\left(\sum_{j=0}^\infty \frac{(e^{-c})^j}{j!}-1\right)=\half\left(e^{e^{-c}}-1\right)\,\,\,\bullet
\end{align*}
\end{proof}


\chapter{The Cut-Off Phenomena}
\section{Introduction}
Given an ergodic random walk $\xi$, a number of techniques for bounding $\|\nu^{\star k}-\pi\|$ have been developed. Recall the mixing time, $\tau$, as the minimum $k$ such that $\|\nu^{\star k}-\pi\|\leq 1/2e$. In particular, as $\|\nu^{\star k}-\pi\|$ is decreasing in $k$, if $\|\nu^{\star k}-\pi\|\leq1/2e$, then $\tau\leq k$. In many random walks, behaviour called the \textit{cut-off phenomenon} occurs and it makes sense to talk about the mixing time, $\tau$, as the time when $\xi$ is random.

\begin{figure}[h]\cone\epsfig{figure=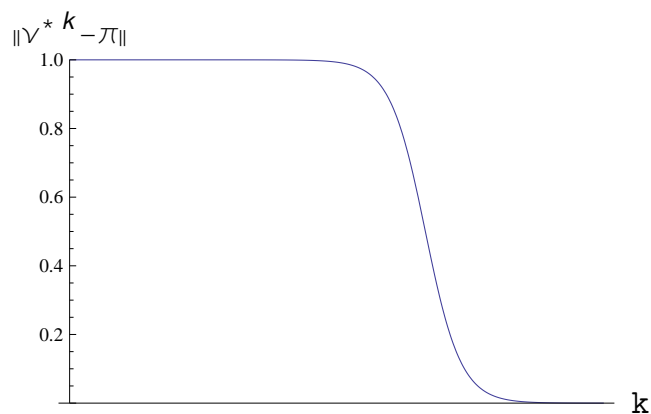}\ctwo\caption{In the cut-off phenomenon, variation distance remains close to 1 initially until the mixing time $\tau$ when it rapidly converges to 0.}\end{figure}

\newpage
In the cut-off phenomenon, the random walk remains far from random until a certain time when there is a phase transition and the random walk rapidly becomes close to random.
\wxn\subsection{Example: Random Transpositions}\wxo
As described in Section \ref{rwfgi}, repeated random transpositions of $n$ cards can be modelled as repeatedly convolving the measure:
\beqa
\nu_n(s):=\left\{\begin{array}{cc}
1/n & \text{ if }s=e
\\ 2/n^2 & \text{ for $s$ a transposition }
\\ 0 & \text{ otherwise }
\end{array}\right.
\enqa
Careful analysis of the representation theory of the symmetric group and an application of the Upper Bound Lemma yields \cite{PD}, for $k=(n\log n)/2+cn$, for $c>0$:
\beq
\|\nu_n^{\star k}-\pi_n\|\leq ae^{-2c}
\enq
for some constant $a$. For a lower bound, Diaconis considers the set $A\subset S_n$ of permutations with one or more fixed points. Two classical results of Feller\footnote{namely the \emph{matching problem}  and the computation of the probability that when $2k$ balls are dropped into $n$ boxes, that one or more of the boxes will be empty \cite{Fell}} give sharp approximations of $\nu_n^{\star k}(A)$ and $\pi_n(A)$ and hence a lower bound for the variation distance may be given. For $k=(n\log n)/2-cn$, $c>0$, as $n\raw\infty$:
\beq
\|\nu_n^{\star k}-\pi_n\|\geq \left(\frac{1}{e}-e^{-e^{-2c}}\right)+o(1)
\enq
Hence for $n$ large, the random walk experiences a phase transition from order to random at $t_n=n \log n/2$. Indeed, this was the first problem where a cut-off was detected (\cite{dash}).

\newpage
\section{Formulation}
There are a number of roughly equivalent formulations of the cut-off phenomenon. The subject developed from the question \emph{how many times must a deck be shuffled until it is close to random?} Card shuffling is modelled by a random walk on $S_n$ where the shuffle is defined by the driving probability $\nu\in M_p(G)$. In most cases, the driving probability $\nu$ is related to $n$ so it makes sense to talk about a natural family of random walks $(S_n,\nu_n)$. When a good asymptote of  the mixing times of these walks was accessible,  it was found that in a number of examples that the cut-off behaviour becomes sharper as $n\raw \infty$. As a corollary of this development, the cut-off phenomenon is defined with respect to the limiting behaviour of a \emph{natural family} $(G_n,\nu_n)$.

\bigskip

 In general, a formulation will be referenced to a particular distance of closeness to random. Surprisingly,  given different norms on $M_p(G)$, a random walk  exhibiting the cut-off phenomenon in the first need not exhibit the cut-off phenomenon in the second. There are a number of roughly equivalent formulations (see Chen's thesis \cite{chen}) that introduce a window size $w_n$. This means that the variation distance goes from 1 to 0 in $w_n$ steps rather than 1 however these formulations still require that $\tau_n\gg w_n$ such that $w_n/\tau_n\raw 0$ hence there is still abrupt convergence.  The original formulation of  Aldous \& Diaconis \cite{chenthree} appeals to an arbitrary sharpness of  convergence of variation distance  to a step function:
\wxn
\subsection{Definition}\wxo
A family of random walks $(G_n,\nu_n)$ exhibits the \textit{cut-off phenomenon} if there exists a sequence of real numbers $\{t_n\}_{n=1}^\infty$ such that given $0<\eps<1$, in the limit as $n\raw\infty,$ the following hold:
\begin{itemize}
\item[(a)]  $\|\nu_n^{\star \lfloor(1+\eps)t_n\rfloor}-\pi_n\|\raw 0$
\item[(b)] $ \|\nu_n^{\star \lfloor(1-\eps)t_n\rfloor}-\pi_n\|\raw1$
\item[(c)] $t_n\raw \infty$
\end{itemize}

If $\tau_n$ is the mixing time of $(G_n,\nu_n)$ presenting cut-off, then the above formulation implies that $\tau_n\sim t_n$ so it makes sense to say that $t_n$ is the time taken to reach random.

%

\subsubsection{Example: Walk on the $n$-Cube}
Recall the walk on the $n$-Cube from the last chapter. Along with the upper bound extracted from the Diaconis-Fourier theory, tedious but elementary calculations bound the variation distance away from 0 for $k=(n+1)(\log n-c)/4$ for $n$ large and $c>0$ (\cite{cecc} --- Th. 2.4.3). This is done via the test function $\phi(s)=n-2w(s)$ whose expectation and variance under $\pi$ are easy to calculate (namely 0 and $n$). The set $A_\beta\subset \Z_2^n$ is essentially defined as the elements whose weight is sufficiently close to $n/2$ for some $\beta$:
\beqa
A_\beta:=\{s\in\Z_2^d:|\phi(s)|<\beta \sqrt{n}\}
\enqa
Use of the Markov inequality bounds $\pi_n(A_\beta)$ above $1-1/\beta^2$. More intricate calculations yield $\nu_n^{\star k}(A_\beta)\leq 4/\beta^2$ and thence
\beq
\|\nu_n^{\star k}-\pi_n\|\geq 1-\frac{5}{\beta^2}
\enq
A more precise definition of $\beta$ in terms of $c$ makes this lower bound useful\footnote{if $\beta=e^{c/2}/2$ then the lower bound is $1-20/e^c$, which clearly tends to $1$ as $c$ increases}. Hence it follows that the random walk has a cut-off at time $t_n=n\log n/4$.
\newpage
\subsubsection{Example: Simple Walk on the Circle}
The simple walk on the circle does not exhibit cut-off. Considering the bounds developed in Section \ref{circ}, note that at $k=n^2/2$, $\|\nu_n^{\star k}-\pi_n\|\leq e^{-\pi^2/4}$, and due to the decreasing nature of $\|\nu_n^{\star k}-\pi_n\|$ this is an upper bound for all $k\geq n^2/2$. Similarly at $k=3n^2/2$:
\beqa
\|\nu_n^{\star k}-\pi_n\|\geq \half e^{-3\pi^2/4-3\pi^4/4n^2}\underset{n\raw \infty}{\raw}\half e^{-3\pi^2/4}
\enqa
and this lower bound holds for all $k\leq 3n^2/2$.

\begin{figure}[h]\cone\epsfig{figure=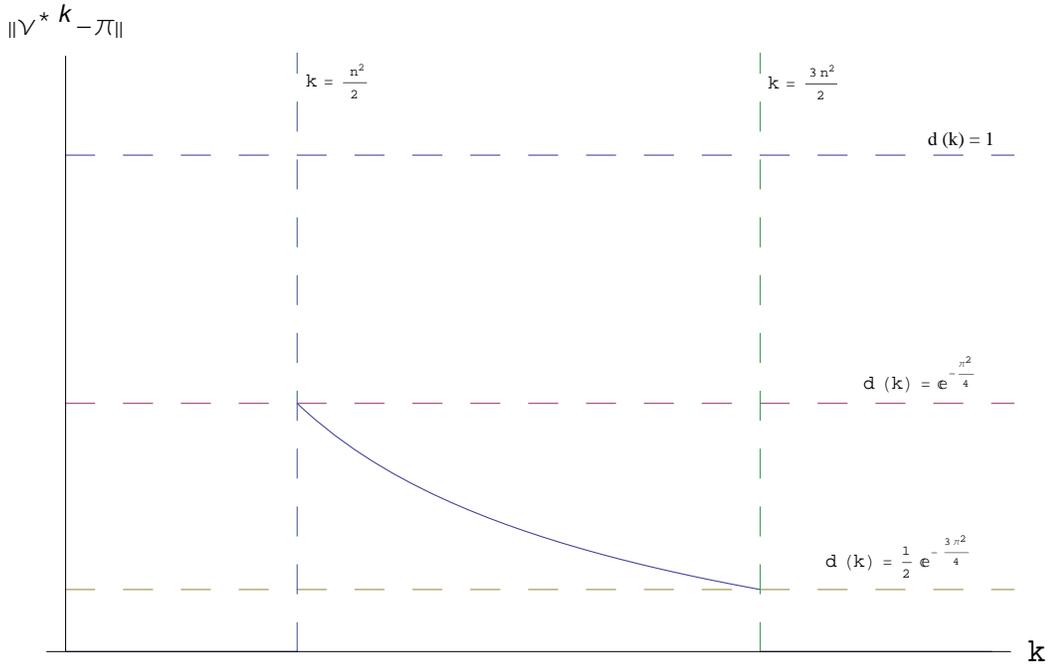}\caption{In the limit as $n\raw \infty$ the simple walk on the circle does not experience an abrupt transition from far from to close to random. Note that $d(k):=\|\nu_n^{\star k}-\pi_n\|$ and the graph is not to scale.}\ctwo\end{figure}

\bigskip

It is an open problem to determine for which families of random walks $(G_n,\nu_n)$ does cut-off occur. Unfortunately there does not appear to be a nice condition for an isolated random walk $\xi$ to exhibit cut-off. In contrast,  given $G$ and $\nu\in M_p(G)$, the ergodic theorem \ref{eth} determines whether or not $(G,\nu)$ is ergodic.


\newpage

An initial attempt at reformulation would be to have as fundamental a period of `far from random' and a period of sharp transition to `close to random'. Rather than being arbitrarily far from random and arbitrarily close to random (in the limit), this finitary formulation would have to define controls $a, b>0$ for far and close to random:
\wxn
\subsection{Definition}
\wxo A random walk on $G$ driven by $\nu\in M_p(G)$ has \textit{$(a,b,q)$ finitary cut-off} if $A:=\{k:\|\nu^{\star k}-\pi\|\geq 1-a\}$, $B:=\{k:b\leq\|\nu^{\star k}-\pi\|\leq 1-a\}$ and $q=|A|/|B|$.

\bigskip

Therefore if $(G_n,\nu_n)$ presents cut-off, each member also has $(a_n,b_n,q_n)$ finitary cut-off, where $a_n,\,b_n\raw0$, $|A_n|\raw\infty$,  and $q_n\raw \infty$. However, consider the natural family $(\Z_n,\nu)$ where $\nu$ is uniform on $\{0,\pm1\}$. This family  has $(1/2,1/4,\mathcal{O}(1))$ finitary cut-off but does not present the cut-off phenomenon. For a family, therefore, presenting cut-off is strictly stronger than presenting finitary cut-off. It is pretty clear that \emph{all} random walks have some level of finitary cut-off. Is there an appropriate level of quality of cut-off?

\begin{figure}[h]\cone\epsfig{figure=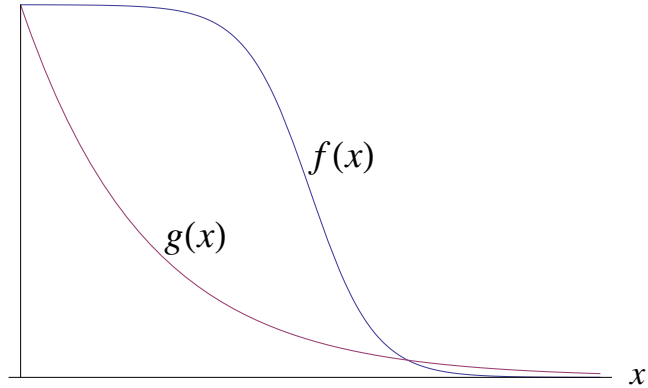}\ctwo\caption{In a natural definition of cut-off, the exponential function $g$ should not have cut-off. The other function, $f$, certainly exhibits some level of cut-off.\label{nat}}\end{figure}

 A continuous version of $(a,b,q)$ finitary cut-off can be considered. Let $f : \R_+\raw [0,1]$ be a non-increasing continuous function with $f(0)=1$ and $f(x)\raw 0$. $f$ exhibits $(a,b,q)$ finitary cut-off where $A=\inf\{x\,:\,f(x)=1-a\}$, $B=\inf\{x\,:\,f(x)=b\}$ and $q=A/(A-B)$.
 \bigskip

In Figure \ref{nat}, $f$ has $(1/2e,1/2e,2.52)$ finitary cut-off while $g$ has $(1/2e,1/2e,0.14)$ finitary cut-off. In a number of examples of established cut-off, e.g. the top-to-random shuffle \cite{chenthirteen}, it has been shown that $\|\nu_n^{\star \lfloor (1-\eps)t_n\rfloor}-\pi_n\|\raw 1$ doubly exponentially as $\eps\raw 1$. Hence consider $(1/e^{2e},1/2e,1)$ finitary cut-off as an appropriate level for cut-off. Indeed  $f$ has $(1/e^{2e},1/2e,0.52)$ finitary cut-off while $g$ has $(1/e^{2e},1/2e,0.0026)$ finitary cut-off. However this too runs into problems. Consider the family of functions $f_d(x)=(1-\tanh(d(x-1/2)))/2$. This family has $(1/e^{2e},1/2e,1)$ finitary cut-off for $d\gtrsim 12.4$.

Diaconis remarks  \cite{PD} that Aldous \& Diaconis  have shown that for most probability measures on a finite group $G$, $\|\nu^{\star 2}-\pi\|\leq 1/|G|$, so for large groups, most random walks are random after two steps.


Therefore, without an alternative formulation of the cut-off phenomenon, it seems likely there will never be a theorem of the form: \emph{A random walk on $G$ with driving probability $\nu\in M_p(G)$ presents `the' cut-off phenomenon at time $k$ if and only if property $P$ is satisfied.}

\section{What Makes it Cut-Off?}
To demonstrate the intransigence of the problem note that the asymptotics of a reversible random walk $\|\nu^{\star k}-\pi\|\sim C\lambda_\star^k$ cannot detect cut-off. A critical idea for understanding of the cut-off phenomena is that variation distance is sensitive. Suppose a deck of cards is shuffled (by $\nu\in M_p(S_{52})$) but the shuffle leaves the ace of spades at the bottom of the deck. If $A\subset S_{52}$ are the arrangements of the deck with the ace of spades at the bottom, then $\nu(A)=1$ but $\pi(A)=1/52$ and $\|\nu-\pi\|\geq 1-1/52$; the deck is very far from random in variation distance! Similarly suppose that  after shuffling by $\nu$ that the ace of spades is in the bottom half of the deck. By letting $B\subset S_{52}$ be all such arrangements it is clear $\|\nu-\pi\|\geq 1/2$. So for any shuffle the entire deck must be well shuffled; it won't do to have even coarse information on  a single card.

\newpage

To illustrate further, consider the top-to-random shuffle. This is the shuffle that takes the top card of the deck and inserts it back into the deck randomly\footnote{i.e. driven by the measure constant on the cycles $(1,m,m-1,\dots,3,2)$, $m=1,\dots,52$}.  Suppose the initial arrangement has the ace of spades  at the bottom of the deck. Initially it will take a while for a card from the top to be placed underneath the ace of spades but eventually one will be and the ace of spades will be second from bottom. After a great number of shuffles the ace of spades will eventually surface at the top of the deck. At every stage up to this point, to within a statistical deviation, the ace of spades is in a specific portion of the deck, dependent on the number of shuffles. Hence up to this point the deck will be far from random. After this step however the ace of spades shall be placed at a random position in the deck and there is every chance the deck is random. It will be seen in the next chapter that the time for the bottom card to come to the top is essentially the time to random and hence the cut-off time.

\bigskip

The survey article by Diaconis \cite{chen11} suggests a number of reasons why cut-off may occur. Diaconis claims that high-multiplicity of second eigenvalue implies cut-off after a remark of Aldous \& Diaconis \cite{chenthree}. The result from \cite{sc}
\beq
\|\nu^{\star k}-\pi\|_2\geq m_\star \lambda_\star
\enq
has some implications for this claim in the two norm (see Chen \cite{chen}). However, in this thesis, cut-offs in  variation norm are the subject of study. One might fear  `folklore heuristic' failure here.   Indeed the claim of Diaconis is almost cited as fact by Hora \cite{jtwo,jone}. Perhaps a more measured statement would be that to show cut-off the random walk may have to exhibit a high degree of symmetry which can imply high multiplicity of the second largest eigenvalue. In the extreme case of almost all eigenvalues equal to $\lambda_\star$ (remembering the average of the eigenvalues is $\nu(e)$), the variation distance looks like $C\lambda_\star^k$ and this doesn't look like cut-off.

\newpage

Chen \cite{chen} discusses a conjecture of Peres that a general Markov chain exhibits the cut-off phenomenon if and only if $\tau_n(1-\lambda_{n\star})\underset{n\raw \infty}{\raw}\infty$. Any Markov chain with cut-off will satisfy this condition. Chen \& Saloff-Coste \cite{csc} have proved this conjecture in the $p$-norm case for $1<p<\infty$ however Aldous has given a Markov chain which is a counterexample in variation distance \cite{chen}. Presently there is no known counterexample to Peres' conjecture in the case of random walks on groups.

\bigskip

Theorem \ref{12}  is relevant for  family of groups $(G_n,\nu_n)$ of moderate growth
with $|\Sigma|$, $A$, $d$ fixed as $n\raw\infty$. These random walks take  large multiple of $\Delta_n^2$ to get random. While a small multiple of $\Delta_n^2$ is not sufficient for randomness, the
transition from 1 to 0 as the number of steps grows  is smooth so that the cut-off is not exhibited \cite{chenseventeen}.  Diaconis \cite{chen11} notes that --- via Gromov's Theorem for  nilpotent groups of finite index --- this result is generic.  For random walks on families of nilpotent groups where $|\Sigma|$ and the index are bounded as $n\raw\infty$,  order $\Delta^2$ steps are necessary for convergence and there is no cut-off. Two examples of such walks are the simple walk on the circle and the walk on the Heisenberg groups, and indeed these are the canonical examples where cut-off does not occur.
\chapter{Probabilistic Methods}

\section{Stopping Times}
In previous chapters the convergence behaviour of a random walks has been examined. It is natural to ask questions of the type \emph{from which time $T$ onwards does $\xi_T$ have a particular property}.  As a simple example of such a random time, consider a random walk $\xi$. The lowest $T_0$ such that $\xi_{T_0}=e$ is such a random time, namely the \textit{first return time}.

To make precise, let $\mathcal{A}_k$ be the $\sigma$-algebra generated by the random variables $\{\xi_j:j\leq k\}$, for   $j,\,k\in\N_0$. Then the $\sigma$-algebra generated by the $\sigma$-algebras $\{\mathcal{A}_k:k\in\N_0\}$, $\mathcal{A}$, canonically admits an increasing sequence:
 $$\mathcal{A}_0\subset \mathcal{A}_1\subset\cdots\subset \mathcal{A}_k\subset\cdots\subset \mathcal{A}$$
 of sub-$\sigma$-algebras of $\mathcal{A}$ (i.e. a \textit{filtration}). If $S(G)$ is the set of sequences in $G$, then a \textit{stopping time} is a map $T:S(G)\raw\N\cup\{\infty\}$ which satisfies $\{T\leq k\}\in\mathcal{A}_k$ for all $k\in\N$.

\bigskip

To formalise the first example of a stopping time, the first return time, write $T_0=\min\{k\geq1\,:\,\xi_k=e\}$.
Of course this generalises easily to another example of a stopping time, namely the \textit{first hitting time}, $T_g=\min\{k\geq0\,:\, \xi_k=g\}$.
  More generally, a subset $A\subset G$ has first hitting time $T_A=\min\{k\geq0\,:\,\xi_k\in A\}$

  \newpage

New stopping times may be constructed from old. If $T$ and $S$ are stopping times for a random walk $\xi$, then so are $\min\{T,S\}$, $\max\{T,S\}$, and $T+n$, $n\in\N$ (see \cite{RS} for proof). The standard analysis of stopping times involves an examination of their expectation, $E_\mu$.  There is a strong relationship between the random distribution $\pi$ and stopping times which is given in the following proposition.
  \wxn
  \subsection{Proposition\label{a1}}  \wxo
  \emph{Let $\xi$ be a random walk on a group $G$. Let $T\in\N$ be a non-zero stopping time such that
  $\xi_T=e$ and $E_\mu T<\infty$. Let $g\in G$. Then
  $$E_\mu(\text{number of visits to $g$ before time $T$})=E_\mu T/|G|$$}
\begin{proof}
Taking the approach of \cite{prel} (Proposition 4, Chapter 2),\new write $\rho(g)=E_\mu(\text{number of visits to }g\text{ before time }T)$. Now $$\lambda(g):=\frac{\rho(g)}{\sum_t\rho(t)}=\frac{\rho(g)}{E_\mu T}$$
is a probability measure on $G$. Next it is claimed that  \beq
\sum_{t\in G}\lambda(t)p(t,g)=\lambda(g).\enq
To see this note that
\begin{align*}
\lambda(g)&=\frac{1}{E_\mu T}\sum_{k=0}^\infty\mu(\xi_k=g,T>k).
\end{align*}
If $g=e$, then $\mu(\xi_0=e)=\mu(\xi_T=e)=1$. Also, for $g\neq e$, by hypothesis, $\mu(\xi_0=g)=\mu(\xi_T=g)=0$. Therefore, in the reindexing $\xi_k\raw\xi_{k+1}$,  the term $\mu(\xi_0=g)$ is replaced by $\mu(\xi_T=g)$ (in the event $T=k+1$). Thus
\begin{align*}
\lambda(g)&=\frac{1}{E_\mu T}\sum_{k=0}^\infty\mu(\xi_{k+1}=g,T>k)
\\&=\frac{1}{E_\mu T}\sum_{k=0}^\infty\sum_{t\in G} \mu(\xi_k=t,T>k,\xi_{k+1}=g)
\end{align*}
By the Markov property,
\begin{align*}
\lambda(g)&=\frac{1}{E_\mu T}\sum_{k=0}^\infty\sum_{t\in G} \mu(\xi_k=t,T>k)p(t,g)
\\&=\sum_{t\in G}\frac{\rho(t)}{E_\mu T}p(t,g)=\sum_{t\in G}\lambda(t)p(t,g)
\end{align*}
Thus it is shown that $\lambda P=\lambda$, and so $\lambda$ is in fact the unique stationary distribution. Consequently
\beqa
\lambda(g)=\frac{\rho(g)}{E_\mu T}=\pi(g) \Raw  \rho(g)=\pi(g)E_\mu T\,\,\,\bullet
\enqa
\end{proof}

\section{Strong Uniform Times}
Consider the following shuffling scheme. Given a deck of $n$ cards in order remove a random card and place it on the top of the deck. Repeat this shuffle until the random time $T$ when every card in the deck has been touched. This $T$ is a stopping time and further every arrangement of the deck is equally likely at this time. Call such a stopping time a \textit{strong uniform time}: a stopping time $T$ such that $\mu(\xi_T=g)=1/|G|$. Diaconis \cite{PD} remarks that this is equivalent to $\mu(\xi_k=g|T\leq k)=1/|G|$.

\bigskip

Aldous \& Diaconis \cite{chentwo} gives a classic account of strong uniform times. For many applications, including the random to top shuffle, the classical coupon collector's problem is required knowledge. Consider a random sample with replacement from a collection of $n$ coupons. Let $T$ be the number of samples required until each coupon has been chosen at least once.
\newpage
\wxn\subsection{Coupon Collector's Bound}\wxo
\emph{In the notation above, let $k=n\log n+cn$, with $c>0$. Then
\beq
\mu(T>k)\leq e^{-c}
\enq}
\begin{proof}
The proof is standard but this is taken from \cite{PD}. For each coupon $b$, let $A_b$ be the event \emph{coupon $b$ is not drawn in the first $k$ draws}. The probability of not picking $b$ once is $1-1/n$, hence $\mu(A_b)=(1-1/n)^k$. Thence
\beqa
\mu(T>k)=\mu\left(\bigcup_b A_b\right)\leq \sum_b \mu(A_b)=n\left(1-\frac{1}{n}\right)^k\leq ne^{-k/n}=e^{-c}\,\,\,\bullet
\enqa
\end{proof}

Recall the separation distance $s(k)$. The separation distance is related to strong uniform times via the following theorem:
\wxn\subsection{Theorem}\wxo
\emph{If $T$ is a strong uniform time for a random walk driven by $\nu\in M_p(G)$,  then for all $k$
\beq
\|\nu^{\star k}-\pi\|\leq s(k)\leq\mu(T>k)\label{sut}
\enq
Conversely there exists a strong uniform time such that the rightmost  inequality holds with equality.}
\begin{proof}
Variation distance is controlled by separation distance so it suffices to prove the rightmost inequality.
Again taking the approach of \cite{PD}, let $k_0$ be the smallest $k$ such that $\mu(T\leq k_0)>0$. The inequality (\ref{sut}) holds if $k=\infty$ and for $k<k_0$. For $k\geq k_0$, $s\in G$:
\begin{align*}
s(k)&\leq 1-|G|\nu^{\star k}(s)\leq 1-|G|\mu(\xi_k=s\,,\,T\leq k)
\\s(k)&\leq1-\underbrace{|G|\mu(\xi_k=s|T\leq k)}_{=1}\cdot\mu(T\leq k)
\\&\leq 1-\mu(T\leq k)=\mu(T>k)
\end{align*}
See \cite{PD} (Theorem 4, Chapter 4C) for the converse result $\bullet$
\end{proof}

\bigskip

This result along with the coupon collector's bound applies immediately to the random to top shuffle. The upper bound proved here is supplemented by the (tricky) second result from \cite{PD} to yield another example of a random walk exhibiting cut-off:
\wxn\subsection{Theorem}\wxo
\emph{For the random to top shuffle, let $k=n\log n+cn$. Then
\begin{eqnarray}
\|\nu^{\star k}-\pi\|\leq e^{-c}\text{ for }c\geq 0,
\\ \|\nu^{\star k}-\pi\|\raw 1\text{ as }n\raw\infty, \text{ for negative }c=c(n)\raw-\infty
\end{eqnarray}}

\section{Coupling}
Coupling is a theoretically stronger method than that of strong uniform times. A \emph{coupling} takes a random walk $\xi$ along with the random walk $\Pi$ (with random distribution) and couples them as a product process $(\xi,\Pi)$. The interpretation being that the two random walks evolve until they are equal, at which time they couple, and thereafter remain equal. More formally a coupling of a random walk $\xi$ (with stochastic operator $P$) takes a `random' operator $\Gamma$ on $M_p(G)\times M_p(G)$ and uses it as an input into $(\xi,\Pi)$ such that the marginal distribution of the first factor is precisely the distribution of $\xi$. The operator must be random in the sense that $\Gamma (\mu,\pi)=(\mu P,\pi)$.  Hence $\Gamma(\nu^{\star k},\pi)= (\nu^{\star k+1},\pi)$. The operator must act on $M_p(G)\times M_p(G)$ in such a way that the $\xi_k$ begin to match up with the $\Pi_k$ until all the elements lie  along the diagonal: $\xi_T=\Pi_T$. That is after $T$ steps the process will have the same distribution as the second process: that is after the stopping time $k=T$ steps the walk will be random. Call such a $T$ a \textit{coupling time}. For appropriate couplings, the coupling time, $T$, may be calculated. To make this argument precise a lemma from \cite{PD} about marginal distributions is required.
\wxn\subsection{Lemma\label{lemma4}}\wxo
\emph{Let $G$ be a finite group. Let $\mu_1$, $\mu_2\in M_p(G)$. Let $\mu\in M_p(G\times G)$ with margins $\mu_1$, $\mu_2$. Let $\Delta =\{(s,s):s\in G\}$ be the diagonal. Then
\beqa
\|\mu_1-\mu_2\|\leq \mu(\Delta^C)
\enqa}
\begin{proof}
Following Diaconis \cite{PD}, let $A\subset G$. Thus
\begin{align*}
|\mu_1(A)-\mu_2(A)|=&|\mu(A\times G)-\mu(G\times A)|
\\ =&\left|\mu((A\times G)\cap \Delta)+\mu((A\times G)\cap \Delta^C)\right.
\\  &\left.-\mu((G\times A)\cap \Delta)-\mu((G\times A)\cap \Delta^C)\right|
\end{align*}
The first and third quantities in the absolute sign are equal. The second and fourth give a difference of two numbers, both   smaller than $\mu(\Delta^C)$ $\bullet$
\end{proof}

\wxn
\subsection{Corollary: Coupling Inequality}\wxo
\emph{If $T$ is a coupling time for a random walk driven by $\nu\in M_p(G)$,  then for all $k$}
\beq
\|\nu^{\star k}-\pi\|\leq \mu(T>k)
\enq
\emph{Conversely there exists a coupling such that the  inequality
holds with equality.}
\begin{proof}
Let $\mu$  be the distribution of $(\xi_k,\Pi)$. Then $\mu$ has marginal distributions $\nu^{\star k}$ and $\pi$. Lemma \ref{lemma4} implies that
\beqa
\|\nu^{\star k}-\pi\|\leq \mu(\Delta^C)=\mu(T>k)
\enqa
See \cite{scon} for a proof and discussion of the existence of an optimal coupling time $\bullet$
\end{proof}

\wxn\subsection{Example: A Walk on the $n$-Cube \cite{coup}}\wxo
Consider the walk on $\Z_2^n$ driven by the measure:
\beq
\nu_n(s):=\left\{\begin{array}{cc}
1/2 & \text{ if }s=e
\\ 1/2n & \text{ if }s=e_i\text{ for some }i
\\ 0 & \text{ otherwise.}
\end{array}\right.
\enq
An equivalent formulation is that a coordinate is chosen independently from $\{1,\dots,n\}$ and a coin flip determines whether the coordinate is flipped or not. Consider the following coupling operator $\Gamma$. Suppose $\xi_k=\sum_i\alpha_ie_i$ and coordinate $j$ is chosen at random. If the coin is heads, then $\xi_{k+1}=\sum_{i\neq j}\alpha_ie_i+(1-\alpha_j)e_j$ and the $j$th coordinate of $\Pi_{k+1}=(1-\alpha_j)$. If the coin is tails, $\xi_{k+1}=\xi_k$ but the $j$th coordinate of $\Pi_{k+1}=\alpha_j$. From the marginal viewpoint of $\xi$, $\Gamma$ is identical to sampling by $\nu_n$. It remains to show that the coupling is suitably random (as described above). Suppose coordinate $j$ is chosen. The distribution of each coordinate of $\Pi_k$ is uniform on $\{0,1\}$. Suppose without loss of generality that  the $j$th coordinate of $\xi_k$   is 1. With equal probability the $j$th coordinate of $\Pi_{k+1}$ will be 0 or 1 by the coin flip, hence the coupling operator is suitably random. Hence the coupling time is when all of the coordinates $\{1,\dots,n\}$ have been chosen. The bound on the coupon collector's bound and the coupling inequality implies the walk is random after $n\log n$ steps.

\chapter{Some New Heuristics}
\section{The Random Walk as a Dynamical System}
Although the dynamics of a particle in a random walk are indeed random, the dynamics of its probability
distribution certainly are not. Indeed note the probability distributions $\{\nu^{\star k}\}_{k\in\N}$ evolve deterministically as $\{\delta^e P^k:k\in \N\}$. Thus the random walk has the structure of a dynamical system  $\{M_p(G),P\}$ with fixed point attractor $\{\pi\}$. The two canonical categories of dynamical systems (for which there is an existing literature of powerful methods e.g. \cite{tt}) are  \textit{topological} and \textit{measure preserving} dynamical systems. Unfortunately at first remove $\{M_p(G),P\}$  appears too coarse and structureless to apply any of these powerful methods. Also the mapping function $P$ is not necessarily invertible and this poses further problems. Indeed in many examples of walks exhibiting cut-off, $P$ may be seen to be singular. Hence the assumption that needs to be made on $P$ to  put a structure on $\{M_p(G),P\}$ sufficient for application of dynamical systems methods to the cut-off phenomenon is overly strict. A more fundamental problem occurs in trying to put the structure of a measure preserving dynamical system on the walk in that if a meaningful\footnote{a measure $\kappa$ wouldn't be very meaningful if $\kappa(M_p(G))=\kappa(\{\pi\})$} measure is put on $M_p(G)$, the fact that $(M_p(G))P^k\underset{k\raw \infty}{\raw} \{\pi\}$ would imply that $P$ is in fact not measure preserving.

\section{Charge Theory\label{charge}}
Two features of the ergodic random walk suggest an obvious generalisation. The first is that a stochastic operator conserves the unit weight of $\mu\in M_p(G)$. Suppose $u\in\R^{|G|}$ is a row vector of weight $q$ in the positive orthant. A normalisation ensures  $u/q\in M_p(G)$ hence $uP/q$ has weight 1 and thus $uP$ has weight $q$. A simple calculation shows that given \emph{any} row vector $u\in \R^{|G|}$ of weight $q$, $uP$ also has weight $q$. Therefore stochastic operators are weight preserving. This immediately implies that the left eigenvectors of an ergodic stochastic operator are of weight zero: $u_i P=\lambda_i u_i$ (except $u_1$ of course).

\bigskip

Secondly an ergodic stochastic operator converges to $U=[1/|G|]$ (the matrix with all entries equal to $1/|G|$), so that given a weight 1 row vector $u$, $uP^n$ converges to $\pi$. In particular, if $\xi_0$ is distributed as any signed probability measure (or \emph{charge}: a signed measure on $G$ such that $\rho(G)=1$) $\rho$, the random walk will still converge to the random distribution. This allows an all manner of generalisations. For example, consider the \textit{signed} stochastic operator $Q=[\rho(ht^{-1})]_{th}$ generated by a signed probability measure $\rho$. Under what conditions will $\delta^e Q^n$ converge to the random distribution?

\section{Invertible Stochastic Operators}
In general a random walk need not start deterministically at $e$, but rather in an initial distribution $\mu=\sum_t \alpha_t\delta^t$. However  $\mu P^n=\sum_t\alpha_t\left(\delta^tP^n\right)$. By right-invariance all the $\delta^tP\raw \pi$ and hence $\mu P^n\raw \pi$ for any initial distribution. In this sense there is a loss of information about initial conditions: the walk forgets where it began, where it was and is totally random. The dynamical systems community make distinctions between the behaviour of invertible and non-invertible maps, however this approach has not been exploited for the case of a random walk on a group.
 \newpage It would be desirable to quantify the `folklore thesis' that \cite{KB}:
\begin{quote}
\emph{The loss of information about initial conditions, as the iteration process proceeds in a chaotic regime, is associated with the \text{non-invertibility} of the mapping function... Hence system memory of initial conditions becomes blurred.}
\end{quote}

Consider the case of a singular and symmetric stochastic operator $P$. The spectral theorem implies $\R^{|G|}$ has a basis of (left) eigenvectors of $P$. Hence $\R^{|G|}$ has an eigenspace decomposition $\bigoplus V_t$, where $V_t:=\ker(\lambda_t I-P)$, where $\{\lambda_t:t\in G\}$ are the eigenvalues of $P$ (with the convention $\lambda_1=1$). Consider $M_p(G)\subset \R^{|G|}=\bigoplus_tV_t$. With a non-trivial kernel $P$, can `\emph{destroy information}' and the na\"{i}ve reaction to this would be to  consider $\nu^{\star k}\in V_1\oplus \text{ker }P$ such that $\|\nu^{\star k}-\pi\|\approx1$. Then $\|\nu^{\star k+1}-\pi\|=0$ and there is cut-off. However given $\delta^e\in \bigoplus V_t$, clearly $P$ kills the $\text{ker }P$ terms at the very first iterate, $\delta^e P$, so this heuristic is incorrect. However in contrived examples the sampling could be done by $\nu_1$ until $\nu_1^{\star k}\in V_1\oplus \text{ker }P_2$ but far from random then sampling by $\nu_2$ (or multiplying by $P_2$) would project onto $V_1$. See Section \ref{securb} for more.
\wxn
\subsection{Proposition\label{inv}}\wxo
\emph{A stochastic operator $P$ is invertible if and only if the equation $uP=\pi$ has the unique solution $u=\pi$.
\new If $P$ is an invertible stochastic operator then the following hold:
\begin{itemize}
\item[(i)] If $u$ is an eigenvector of $P$, then $u$ is an eigenvector of $P^{-1}$. In particular, $\pi P^{-1}=\pi$ and $P^{-1}k=k$ for any constant function $k\in F(G)$.
\item[(ii)] If $\{\lambda_t:t\in G\}$ are the eigenvalues of $P$, then $\{1/\lambda_t:t\in G\}$ are the eigenvalues of $P^{-1}$. In particular, 1 is an eigenvalue of $P^{-1}$, and all other eigenvalues of $P^{-1}$ have modulus greater than 1.
\item[(iii)] The signed probability measures on $G$, $M_1(G)$, are stable under $P^{-1}$.
\item[(iv)] For $k\in\N$, $\delta^e P^{-k}\in M_1(G)\bs M_p(G)$.
\end{itemize}}
\begin{proof}
If $P$ is invertible $uP=\pi$ has unique solution. If $P$ is singular then the kernel is non-trivial. Let $u_1\neq u_2\in\text{ker }P$ be normalised such that $\nu_i:=\pi+u_i\in M_p(G)$, then $\nu_iP=\pi$.

\bigskip

(i) and (ii)  are basic linear algebra facts.
\begin{itemize}
\item[(iii)] From (i)  the row and column sums of $P^{-1}$ are 1. Thence let $v\in M_1(G)$;
\begin{align*}
vP^{-1}(G)&=\sum_{s\in G}\left(\sum_{t\in G}v(t)p^{-1}(t,s)\right)
\\&=\sum_{t\in G}v(t)\underbrace{\left(\sum_{s\in G}p^{-1}(t,s)\right)}_{=1}\,\,\,\bullet
\end{align*}
\item[(iv)] From (iii), $\delta^eP^{-1}\in M_1(G)$. Assume there exists $\nu\in M_p(G)$ such that $\nu P=\delta^e$. Now $\nu P(s)=\la \nu,p_s\ra$ must equal $\delta^e(s)$ where $p_s$ is the row vector equal to the $s$-column of $P$. By Cauchy-Schwarz:
\beq
|\la \nu,p_s\ra|\leq \|\nu\|_2\|p_s\|_2\leq\|\nu\|_1\|p_s\|_1
\enq
Because
\beqa
\nu P(e)=\la \nu,p_e\ra=1=\|\nu\|_1\|p_e\|_1
\enqa
the second and third inequalities are equalities for $s=e$. The first equality implies that $\nu$ and $p_e$ are linearly dependent, $\nu=k p_e$. As probability measures must have weight 1, this implies $\nu=p_e$. The second equality implies that $\nu$ and $p_e$  are Dirac measures. Hence $\nu$ is a Dirac measure, say $\delta^g$, and thus $P$ is not ergodic (as $\Sigma$ is a subset of the coset $\{e\}g$,  of the proper normal subgroup $\{e\}$). Inductively given $v\in M_1(G)\bs M_p(G)$, there does not exist $\nu\in M_p(G)$ such that $\nu P=v$ as $v$ must have negative entries but both $\nu$ and $P$ are positive $\bullet$

\end{itemize}
\end{proof}
\newpage
\section{Convolution Factorisations of $\pi$\label{securb}}
Take a deck of cards and transpose the top card with a random card. Next transpose the second card with a random card (at or underneath the second) and continue inductively until all but the second from bottom card has been transposed. Apply the same shuffle to the 51st card ((51,51) or (51,52)). The first card is random, the second is random and inductively all the cards are random. Hence considering the group $S_n$ and the measures $\nu_i$ uniform on the transpositions $\{(i,i),(i,i+1),\dots,(i,n)\}$ the random distribution factorises as:
 \beq
\pi=\nu_{n-1}\star\cdots\star\nu_2\star\nu_{1} \label{urban}
\enq
Urban \cite{urban} considers the question: given a group $G$ and a symmetric set of generators $\Sigma$, does there exist a finite number of convolutions of symmetric measures $\{\nu_i\in M_p(G):i=1,\dots,m\}$ supported on $\Sigma$ such that (\ref{urban}) holds (with $m$ rather than $n-1$ terms)?   Urban uses Diaconis-Fourier theory (particularly Lemma \ref{l2}) to show that if, at a non-trivial irreducible representation of $G$, $\rho$, the Fourier transform of $\nu_m\star\cdots\star\nu_1$ is non-zero then (\ref{urban}) cannot hold. Briefly, Lemma \ref{l2} states that at any non-trivial irreducible representation, $\widehat{\pi}(\rho)=0$; and the Fourier transform of $\nu_m\star\cdots\star\nu_1$ is easily computed via the convolution theorem.

\bigskip

If  $\nu^{\star k}=\pi$ for some finite $k\in\N$ then the results of Section \ref{spectral} shows that  $\nu=\pi$.  In particular, as $\nu$ is symmetric, $P$ has an eigenbasis, and $1$ is an eigenvalue of $P$ with multiplicity 1. Suppose for contradiction that $\nu^{\star k}=\pi$ for some $k\in\N$, but $\nu\neq\pi$. Suppose $\delta^e\in V_1\oplus \ker P$; then $\delta^eP=\pi$. However $\delta^eP=\nu\star\delta^e$, however $\nu\star\delta^e=\nu$ and thus $\nu=\pi$. Hence at least one of the eigenvectors in the eigenbasis expansion of $\delta^e$ is associated with a non-zero eigenvalue. Thus hence $\nu^{\star k}\neq \pi$ for any $k\in\N$. Note that each of the $\nu_i$ induces a stochastic operator $P_i$ and (\ref{urban}) is equivalent to
\beq
U=P_mP_{m-1}\cdots P_2P_1\label{urban1}
\enq
Note that $U$ is singular. If each of the  $P_i$ are invertible then so is $U$, a contradiction. Therefore (\ref{urban1}) cannot be true if each of the $P_i$ are invertible. Theorem 6 on page 49 of Diaconis \cite{PD} implies that each eigenvalue of $\widehat{\nu}(\rho)$, where $\rho$ is an irreducible representation, is an eigenvalue of multiplicity $d_\rho$. In the case of an Abelian group, the eigenvalues of $P$ are simply given by $\{\widehat{\nu}(\rho_i):\rho_i\text{ irreducible}\}$ and the analysis breaks down to that of Urban's as  $\widehat{\nu}(\rho_i)\neq 0$  is equivalent to $0$ is not an eigenvalue of $P$; i.e. $P$ is invertible.
\subsubsection{Example: Simple Walk on the Circle}
Let $n$ be odd and consider the set $\mathcal{M}$ of not-necessarily symmetric measures  with support $\Sigma=\{\pm 1\}$ (i.e. $\mathcal{M}=\{\nu_p\in M_p(G):\nu_p(1)=p,\,\nu_p(-1)=1-p;p\in(0,1)\}$). Does $\pi$ admit a finite convolution factorisation of measures from $\mathcal{M}$? For convenience denote $q:=1-p$ and $\alpha:=p/q$. Consider the  stochastic operator associated to $\nu_p$:
\beqa
P_p=\left(\begin{array}{ccccccc}
0 & p & 0 & 0 & \cdots & 0 & q
\\q & 0 & p & 0 & \cdots & 0 & 0
\\0 & q & 0 & p & \cdots &0 &0
\\ \vdots & & \vdots & & \ddots & & \vdots
\\ p & 0&  0 & 0 & \cdots & q & 0
\end{array}\right)
\enqa
Apply the  elementary row operation $r_i\raw r_i/q$ to each row and permute the rows by $(r_nr_{n-1}r_{n-2}\cdots r_1)$:
\beqa
P_p\equiv\left(\begin{array}{ccccccc}
1 & 0 & \alpha & 0 & \cdots & 0 & 0
\\0 & 1 & 0 & \alpha &  & 0 & 0
\\ &  & \vdots &  &  & &\vdots
\\ \alpha & 0 & 0 & 0 & & 1 & 0
\\ 0 & \alpha &  0 & 0 & \cdots & 0 & 1
\end{array}\right)
\enqa
Now\footnote{if $p<q$ then $\alpha<1$ and Gershgorin's Theorem implies that $P_p$ is invertible. If $p>q$, then $\alpha>1$ and elementary row operations give $P_p$ invertible similarly. Gershgorin cannot deal with the case $p=q$ however. Gershgorin can show $P_p$ is invertible with $n$ even when $p\neq q$, but on this support, the walk is not ergodic.} eliminate by $r_{n-1}\raw r_{n-1} -\alpha r_1$ and $r_{n}\raw r_n-\alpha r_2$:
\beqa
P_p\equiv\left(\begin{array}{ccccccc}
1 & 0 & \alpha & 0 & \cdots & 0 & 0
\\0 & 1 & 0 & \alpha &  & 0 & 0
\\\vdots &  &  &  &  & &\vdots
\\ 0 & 0 & -\alpha^2 & 0 & & 1 & 0
\\ 0 & 0 &  0 & -\alpha^2 &  & 0 & 1
\end{array}\right)
\enqa
Now suppose $n=2m+1$ and continue inductively until:
\beqa
P_p\equiv\left(\begin{array}{ccccc}
1 & & 0 & 0 & 0

\\ &  &  &   &\vdots
\\ 0 &  & (-1)^{m+1}\alpha^m & 1 & 0
\\ 0 & \cdots &  0 & (-1)^{m+1}\alpha^m & 1
\end{array}\right)
\enqa
A final application of $r_{n-1}\raw r_{n-1}-(-1)^{m+1}\alpha^m r_{n-2}$ and $r_{n}\raw r_n-(-1)^{m+1}\alpha^m r_{n-1}$ yields:
\beqa
P_p\equiv\left(\begin{array}{cccc}
1 & & 0 & 0
\\ & \ddots &  &
\\0 &  & 1 & (-1)^{m+2}\alpha^{m+1}
\\ 0 & \cdots & 0 & 1
\end{array}\right)
\enqa
Hence the $P_p$ have $n$ pivots and are thus invertible so a finite convolution of measures from $\mathcal{M}$ is never random.

\bigskip

Urban proves a stronger result using the Diaconis-Fourier theory; namely if $\mathcal{M}$ is a set of measures  symmetric on $\{s\in\Z_n:|s|<n/4\}$ then there is no $\pi$-factorisation. A quick look at the representation theory of $\Z_n$ shows that the Fourier transform of these measures is bounded away from 0 and hence so are the eigenvalues.
\subsubsection{Example: Urban's Transposition Shuffle}
Consider the convolution described by at the start of this section. The final  driving measure $\nu_{n-1}=(\delta^e+\delta^{(n-1,n)})/2$ generates a singular stochastic operator $P_{n-1}$ by Proposition \ref{inv} (v) and a slight rearrangement shows that all of the $\nu_i$ generate singular stochastic operators.
\subsubsection{Open Problem}
This leads onto the interesting  question:
\begin{quote}
\emph{For what measures $\nu\in M_p(G)$ is the associated stochastic operator invertible?}
\end{quote}
A sufficient condition for invertibility guaranteed by Gershgorin's circle theorem is that $\nu(e)>1/2$.

\bigskip

\newpage
\section{Geometry of the $\|\nu^{\star k}-\pi\|$ Graph}
Consider an invertible symmetric ergodic stochastic operator $P$. Due to the fact that the  eigenvalues of $P^{-1}$ (except 1) are all modulus greater than 1, the sequence  $\|\nu^{\star (-k)}-\pi\|$ is monotonically increasing to infinity as $k\raw \infty$. Hence the graph looks something like:

\begin{figure}[h]\cone\epsfig{figure=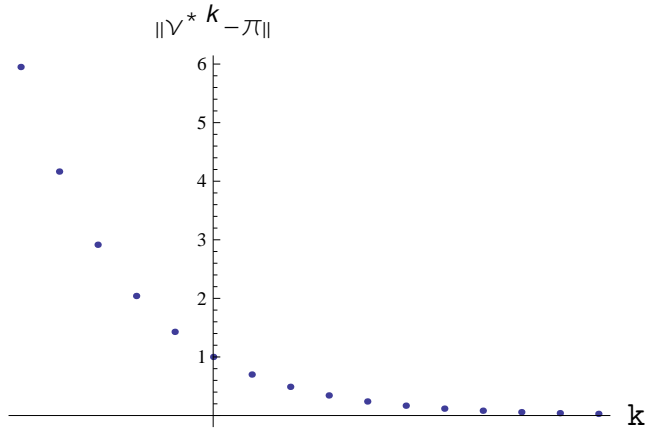}\caption{As $k\raw-\infty$, $\nu^{\star k}$ leaves $M_p(G)$ and becomes a `big' signed measure. }\ctwo\end{figure}

The assumption could be made that in this case the graph \emph{must} be `concave up' and similarly to $g(x)$ in Figure \ref{nat}, does \emph{not} exhibit cut-off. Suppose an invertible stochastic operator did show cut-off:

\begin{figure}[h]\cone\epsfig{figure=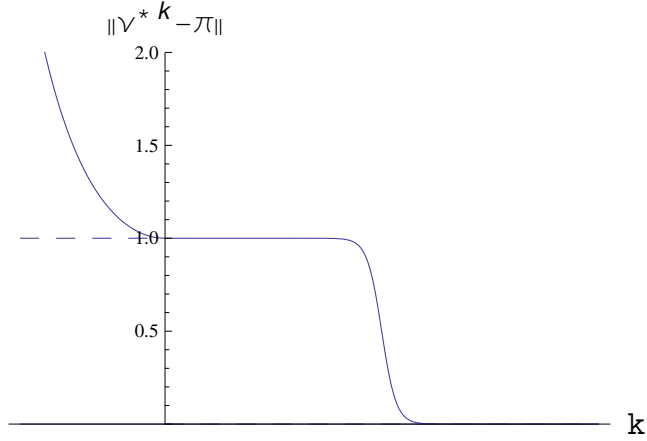}\caption{One could conjecture that the non-dashed line behaviour, supposedly corresponding to an invertible stochastic operator, with two `turning points' be impossible. }\ctwo\end{figure}

Instead one might think that somehow the dashed line behaviour is necessary for cut-off to hold --- and of course this behaviour cannot hold when $P$ is invertible. This leads to the conjecture: $P$ invertible implies no cut-off.  However, in general, $P^{-m}(\{\delta^e\})$ is non-empty, and if a representative $u_m$ from this set is chosen the graph of $\|u_mP^{m+k}-\pi\|$ will exhibit the `non-dashed line' behaviour. Note that for the random walk on the cube with loops there is no charge that is sent to $\delta^e$ by $P$.  This leads onto another interesting  question:
\subsubsection{Open Problem}
\begin{quote}
\emph{For what singular stochastic operators  $P$ generated by $\nu\in M_p(G)$  does there exist a charge $u\in M_1(G)$ such that $u P=\delta^e$?}
\end{quote}

\bigskip

Unfortunately the stochastic operator for the simple walk with loops on $\Z_2^n$ with $n$ even is invertible. If true the conjecture would have placed the problem in a very precarious position. Suppose  $(G_n,\nu_n)$ is a family exhibiting the cut-off phenomenon (so that the stochastic operator is singular), such that $e\in\text{supp}(\nu_n)=\Sigma_n$. Let $\eps\in(0,1/2)$,  and transform the $\nu_n$ as:
\beq
\nu_n^\prime(s):=\left\{\begin{array}{cc}
\half+\eps & \text{ if }s=e
\\[1ex] \frac{1/2-\eps}{|\Sigma_n|-1} & \text{ if }s\in\Sigma_n\bs\{e\}
\end{array}\right.
\enq
Then by Gershgorin's circle theorem $P^\prime$ would be invertible and hence two random walks with the same support need not exhibit the same behavior: the condition for cut-off to hold would not be on the support only. Unfortunately for those active in the field one would assume the condition is indeed this complex.
\chapter{Appendix}
\wxn\section{Proof of Lemma \ref{zpl}}\wxo
\begin{enumerate}
\item
\textbf{Claim:}
\beq
\left|\cos\left(j\frac{\pi}{n}\right)\right|=\left|\cos\left(l\frac{\pi}{n}\right)\right|\,\,\text{ for any }l\in [j]_n
\enq
Suppose $j\equiv l\text{ mod }n$, where $l\in\{0,1,\dots,n-1\}$, so that $j=l+mn$ for some $m\in\Z$. Then
\begin{align*}
\cos\left(j\frac{\pi}{n}\right)&=\cos\left((l+mn)\frac{\pi}{n}\right)
\\&=\cos\left(\frac{l\pi}{n}+m\pi\right)
\\&=\cos\frac{l\pi}{n}\,\cos m\pi-\sin\frac{l\pi}{n}\,\underbrace{\sin m\pi}_{=0}
\\&=(-1)^m\cos\frac{l\pi}{n}
\end{align*}

\bigskip

Now let $a_t=\cos (\pi t/n)$ and $b_t=\cos (2\pi t/n)$, and note that for $t=1,2,\dots,(n-1)/2$:
\beq
|a_t|=\left\{\begin{array}{cccc}\left|b_{(n+t)/2}\right| & \text{ and } & \left|b_{(n-t)/2}\right|  & \text{ if }t\text{ odd}
\\[2ex] \left|b_{t/2}\right| & \text{ and } & \left|b_{n-t/2}\right| & \text{ if }t\text{ even}
\end{array}\right.
\enq
Hence as $(x)^2=|x|^2$:
\beqa
\sum_{t=1}^{n-1} \cos^{2k}(2\pi t/n)=2\sum_{t=1}^{(n-1)/2}\cos^{2k}(\pi t/n)\,\,\,\bullet
\enqa

\item Let $h(x)=\log\left(e^{x^2/2}\cos x\right)$; so that $h'(x)=x-\tan x$ and $h''(x)=1-\sec^2x$. Thus $h''(x)\leq 0$ on $[0,\pi/2]$ and so with $h'(0)=0$, $h(x)$ is a decreasing function in $x$. In particular, $h(x)\leq h(0)=0$ and as $\log$ is an increasing function, $e^{x^2/2}\cos x\leq 1$, for $x\in[0,\pi/2]$ $\bullet$

\item In the first instance:
\beqa
\sum_{j=0}^\infty e^{-3jx}=\frac{1}{1-e^{-3x}}
\enqa
is a convergent geometric series when $x>0$. Now
\beqa
\sum_{j=0}^\infty e^{-3jx}=\sum_{j=1}^\infty e^{-3(j-1)x}.
\enqa
Also $j^2-1\geq 3(j-1)$ for each $j\in\N_0$.
Hence, as $e^x$ is increasing, for all $j\in\N_0$, $e^{-(j^2-1)x}\leq e^{-3(j-1)x}$, and so
\beqa
 \sum_{j=1}^\infty e^{-(j^2-1)x}\leq \sum_{j=1}^\infty e^{-3(j-1)x}= \sum_{j=0}^\infty e^{-3jx}\,\,\,\bullet\enqa
\item Taking the approach of \cite{cecc}, let $h(x)=\log\left(e^{x^2/2+x^4/2}\cos x\right)$;
\begin{align*}
h(0)&=0
\\  h'(x)&=\left.x+x^3-\tan x\right|_{x=0}=0
\\ h''(x)&=\left.3x^2-\tan^2x\right|_{x=0}=0
\\ h'''(x)&=\left.6x-2\sec^2 x\tan x\right|_{x=0}=0
\\ h^{iv}(x)&=6+4\sec^2x-6\sec^4x
\end{align*}
This is a quadratic in $\sec^2 x$ which is positive when $|\sec x|\leq \sqrt{1+\sqrt{10}/3}$. This translates into better than $x\in[0,\pi/6]$ $\bullet$
\end{enumerate}

\newpage
\wxn\section{Proof of Lemma \ref{hcl}}\wxo
\begin{enumerate}
\item[1.] In the first instance:
\beqa
1-\frac{2(n+1-l)}{n+1}=1-2+\frac{2l}{n+1}=-\left(1-\frac{2l}{n+1}\right)
\enqa
So that
\beqa
\left(1-\frac{2l}{n+1}\right)^{2k}=\left(1-\frac{2(n+1-l)}{n+1}\right)^{2k}
\enqa
Secondly,
\begin{align*}
{n\choose l}-{n\choose n+1-l}&=\frac{n!}{l!(n-l!)}-\frac{n!}{(n+1-l)!(\cancel{n}-(\cancel{n}+1-l))!}
\\ &=\frac{n!}{(l-1)!(n-l)!}\left[\frac{1}{l}-\frac{1}{n+1-l}\right]
\\ &=\frac{n!}{(l-1)!(n-l)!}\left[\frac{n+1-l-l}{l(n+1-l)}\right]
\\ &=\frac{n!}{(l-1)!(n-l)!}\left[\frac{n+1-2l}{l(n+1-l)}\right]\underset{l\leq n/2}{\geq}0
\end{align*}
That is, if $l\leq n/2$,
\beqa
{n\choose l}\geq{n\choose n+1-l}\,\,\,\bullet
\enqa

\item[2.]
By definition,
\beqa
{a\choose b}=\frac{a!}{b!(a-b)!}=\frac{a(a-1)(a-2)\cdots(a-b+1)}{b!}\leq \frac{a^b}{b!}\,\,\,\bullet
\enqa

\item[3.]  It suffices to show
\beq
f(j):=\log\left(1-\frac{2j}{n+1}\right)^{2k}\leq-j\log n-jc=:g(j)
\enq
as $\exp$ is an increasing function. Now writing $k=(n+1)(\log n+c)$,
\begin{align*}
f(1)&=\half(n+1)(\log n+c)\log \left(1-\frac{2j}{n+1}\right)\text{ , and }
\\ g(1)&=-(c+\log n)
\end{align*}
Now $c=4k/(n+1)-\log n$ so $c+\log n=4k/(n+1)$. Therefore
\beqa
f(1)-g(1)=\left(\frac{4k}{n+1}\right)\left[1+\half(n+1)\log \left(\frac{n-1}{n+1}\right)\right]
\enqa
This is negative ($f(1)\leq g(1)$) if
\begin{align*}
1+\half(n+1)\log\left(\frac{n-1}{n+1}\right)&\leq 0
\\\Leftrightarrow \log\frac{n-1}{n+1}&\leq-\frac{2}{n+1}
\\\Leftrightarrow h(n)&=\log \left(\frac{n+1}{n-1}\right)- \frac{2}{n+1}\geq 0
\end{align*}
Now $h(2)=\log 3-1>0$ and
\beqa
\lim_{n\raw\infty}\left[\log\underbrace{\left(\frac{n+1}{n-1}\right)}_{\raw 1}-\frac{2}{n+1}\right]=0.
\enqa
Differentiating with respect to $n$,
\beqa
\\h'(n)=-\frac{2}{n^2-1}+\frac{2}{(1+n)^2}=-\frac{4}{(n+1)^2(n-1)}\leq 0.
\enqa
Hence $h(n)$ is monotone decreasing from $h(2)>0$ to $0$ so is positive. Hence $f(1)\leq g(1)$. Now differentiating with respect to $j$,
\beqa
f'(j)=-\frac{(n+1)(c+\log n)}{n+1-2j}=-\frac{4k}{n+1-2j}\underset{j\leq n/2}{\leq} 0
\enqa
Also
\beqa
\\g'(j)=-c-\log n=-\frac{4k}{n+1}
\enqa
Finally as $j\geq0$, $f'(j)\leq g'(j)$, for all $j\leq n/2$ $\bullet$

\end{enumerate}

\bibliography{mybib}{}
\bibliographystyle{plain}

\end{document}